\documentclass[journal,twoside,web]{ieeecolor}

\usepackage{generic}
\usepackage{mathtools}
\mathtoolsset{showonlyrefs}

\usepackage{cite}
\usepackage{amsmath,amssymb,amsfonts}
\usepackage{algorithmic}
\usepackage{graphicx}
\usepackage{textcomp}
\usepackage[T1]{fontenc}

\usepackage{euscript}
\usepackage{xcolor}
\usepackage{epsfig}
\usepackage{listings}
\usepackage{longtable}
\usepackage{multirow}
\usepackage{epstopdf}
\usepackage{bm}
\usepackage{makecell}
\usepackage{subcaption}
\usepackage{matlab-prettifier}
\usepackage{mathrsfs}

\newtheorem{theorem}{\bf \text{Theorem}}

\newtheorem{assumption}{\bf \text{Assumption}}

\newtheorem{corollary}{\bf \text{Corollary}}[theorem]
\newtheorem{lemma}{\bf \text{Lemma}}
\newtheorem{remark}{\bf \text{Remark}}

\newtheorem{example}{\bf \text{Example}}

\newcommand{\E}{{\mathbb E}}
\newcommand{\COV}{{\mathbb C\mathbb O\mathbb V}}

\newcommand{\R}{{\mathbb R}}
\newcommand{\C}{{\mathbb C}}

\newcommand{\SURE}{\text{SURE}}

\makeatletter

\newcommand{\Rmnum}[1]{\expandafter \@slowromancap\romannumeral #1@}
\makeatother

\DeclareMathOperator*{\rank}{rank}
\DeclareMathOperator*{\diag}{diag}
\DeclareMathOperator*{\cond}{cond}
\DeclareMathOperator*{\Tr}{Tr}

\DeclareMathOperator*{\Var}{Var}
\DeclareMathOperator*{\GCV}{GCV}

\DeclareMathOperator*{\EB}{EB}
\DeclareMathOperator*{\Sy}{Sy}
\DeclareMathOperator*{\tb}{b}

\DeclareMathOperator*{\ty}{y}

\DeclareMathOperator*{\XMSE}{XMSE}
\DeclareMathOperator*{\ML}{ML}

\DeclareMathOperator*{\MSE}{MSE}

\DeclareMathOperator*{\XBias}{XBias}
\DeclareMathOperator*{\XVar}{XVar}
\DeclareMathOperator*{\XVarHPE}{XVarHPE}
\DeclareMathOperator*{\VarHPE}{VarHPE}

\DeclareMathOperator*{\HOT}{HOT}
\DeclareMathOperator*{\Bias}{Bias}
\DeclareMathOperator*{\Bayes}{Bayes}

\DeclareMathOperator*{\TR}{R}
\DeclareMathOperator*{\FIT}{FIT}

\DeclareMathOperator*\argmin{arg\,min}

\DeclareRobustCommand*\C{C}

\allowdisplaybreaks[4]

\begin{document}
\title{Excess Mean Squared Error of Empirical Bayes Estimators}
\author{Yue Ju, Bo Wahlberg, \IEEEmembership{Fellow, IEEE}, and H\r{a}kan Hjalmarsson, \IEEEmembership{Fellow, IEEE}
\thanks{This work was supported by VINNOVA Competence Center AdBIOPRO, contract [2016-05181] and by the Swedish Research Council through the research environment NewLEADS(New Directions in Learning Dynamical Systems), contract [2016-06079], and contract [2019-04956].}
\thanks{Yue Ju, Bo Wahlberg, and H\r{a}kan Hjalmarsson are now with Division of Decision and Control Systems, School of Electrical Engineering and Computer Science, KTH Royal Institute of Technology, 10044 Stockholm, Sweden (e-mails: yuej@kth.se, bo@kth.se, hjalmars@kth.se). H\r{a}kan Hjalmarsson is also with the Competence Centre for Advanced BioProduction by Continuous Processing, AdBIOPRO.}}

\maketitle

\begin{abstract} 
Empirical Bayes estimators are based on minimizing the average risk with the hyper-parameters in the weighting function being estimated from observed data. The performance of an empirical Bayes estimator is typically evaluated by its mean squared error (MSE). However, the explicit expression for its MSE is generally unavailable for finite sample sizes. To address this issue, we define a high-order analytical criterion: the excess MSE. It quantifies the performance difference between the maximum likelihood and empirical Bayes estimators. An explicit expression for the excess MSE of an empirical Bayes estimator employing a general data-dependent hyper-parameter estimator is derived. As specific instances, we provide excess MSE expressions for kernel-based regularized estimators using the scaled empirical Bayes, Stein’s unbiased risk estimation, and generalized cross-validation hyper-parameter estimators. Moreover, we propose a modification to the excess MSE expressions for regularized estimators for moderate sample sizes and show its improvement on accuracy in numerical simulations.
\end{abstract}

\begin{IEEEkeywords}
Empirical Bayes estimators, High-order asymptotic theory, Mean squared error, Hyper-parameter estimator. 
\end{IEEEkeywords}


\section{Introduction}

\IEEEPARstart{B}{ayes} estimation is concerned with estimating parameters that minimize the Bayes risk, which is the weighted average of a given loss function by incorporating the weighting function of parameters and the likelihood function of observed data, see, e.g., \cite{Lehmann:06, B2013}. {In most cases, the finite-sample statistical properties of Bayes estimators are analytically intractable}, while asymptotic theory \cite{Vaart1998, LY2000, L12} provides a powerful framework for the analysis of their large-sample properties.

{For a constant weighting function}, the Bayes estimator is equivalent to the maximum likelihood (ML) estimator. Its consistency and asymptotic normality {have} been discussed under a variety of conditions, e.g., \cite{L1953, C1954, D1961, H1971, Ljung1999, IR2013}. For a general given weighting function, the asymptotic equivalence between the Bayes and ML estimators can be established using the Bernstein-von Mises Theorem \cite[Chapter 10.2]{Vaart1998}. The consistency and asymptotic normality of Bayes estimators for independent identically distributed and independent nonhomogeneous observations have been studied in \cite[Chapter 3.2-3.4]{IR2013}. These asymptotic properties have also been derived for the ergodic diffusion processes \cite{Y2011, O2019}.

The family of Bayes estimators is a subset of the family of empirical Bayes (EB) estimators, whose weighting functions are allowed to be functions of observed data. EB estimation was first introduced in \cite{R1956}. There is a well-established asymptotic theory for this family of estimators, e.g., the convergence of the Bayes risks of EB estimators is studied in \cite{R1964, DZ1976, S1976, WZ95, M2018}. Note that various regularized estimators can be interpreted as EB estimators with different forms of weighting functions. These regularized estimation methods have been increasingly recognized as a complement to classical system identification \cite{PDCDL14, PCCDL22}. 

In particular, for the kernel-based Gaussian weighting function, the corresponding EB estimator is called the generalized Tikhonov regularized estimator. Its asymptotic properties have been discussed in \cite{GS2006, PC15, MCL18asy, MC23, JCML21cdc, JCML22ccc, ZCM24}. For the matrix-variate Gaussian weighting function, the corresponding EB estimator becomes the nuclear norm regularized estimator, and its asymptotic properties have been studied in \cite{MW2018}. For the Laplace weighting function, the corresponding EB estimator is known as the least absolute shrinkage and selection operator (LASSO) regularized estimator, and its asymptotic properties have been studied in \cite{KF2000, BG2010}. For the mixture of the Gaussian and Laplace weighting functions, the corresponding EB estimator is referred to as the elastic net regularized estimator. Its asymptotic properties have been investigated in \cite{UV2010}. 

In this paper, we consider linear regression models and EB estimators with a quadratic loss function. Under certain conditions, the ML and EB estimators can both be proved to be consistent and asymptotically efficient \cite{Ljung1999, JCML21eb}. This means that they have the same first-order asymptotic statistical properties. To distinguish the performance of these estimators, we will utilize high-order asymptotic theory \cite{G94Higher, BGV1997}. 

For the ML estimator, its second-order expansion has been discussed in \cite{Rao61, GS74, G94Higher}. It suggests that {the ML estimator can be improved} by correcting the bias and minimizing the variance introduced by the second-order expansion. The third-order expansion of the ML estimator has been studied in \cite{AT2012}.

For EB estimators, second-order and third-order asymptotic efficiency of generalized Bayes estimators are considered in \cite{TA1979}. For a Gaussian weighting function with zero mean and a covariance matrix parameterized by hyper-parameters, the corresponding EB estimator is referred to as the kernel-based regularized estimator in this paper. Its high-order asymptotic properties depend on the asymptotic properties of its hyper-parameter estimators. The consistency of a general hyper-parameter estimator has been derived in \cite[Section B.2]{MCL18asy} by using \cite[Theorem 8.2]{Ljung1999}. More specifically, the consistency of the EB, the Stein's unbiased risk estimation (SURE), and the generalized cross-validation (GCV) hyper-parameter estimators have been established in \cite{PC15, MCL18asy, MC23}. Moreover, the asymptotic normality of these hyper-parameter estimators has been proven in \cite{JCML21eb, JCML21cdc, JCML22ccc}. Based on these convergence results, the bias and variance introduced by the second-order and third-order expansions of the EB-based, SURE-based and GCV-based regularized estimators for stochastic filtered white noise inputs are further studied in \cite{JCML21eb}.    

In this paper, we consider finite impulse response (FIR) models with deterministic inputs and white Gaussian measurement noise. To analyze the difference in performance of the EB and ML estimators, we introduce the excess mean squared error (XMSE) of an EB estimator. It is defined as the high-order limit of the mean squared error (MSE) difference between the EB and ML estimators. The XMSE is positive when the ML estimator performs better than the EB one for large sample sizes. The main contributions of this work are:
\begin{enumerate}
\item We derive an explicit expression for the XMSE of an EB estimator equipped with a general data-dependent hyper-parameter estimator.
\item As specific instances, we present XMSE expressions for generalized Bayes estimators, and kernel-based regularized estimators with the most common hyper-parameter estimators, i.e., the scaled EB, SURE, and GCV ones.
\item We conduct numerical simulations to demonstrate the accuracy of XMSE expressions for regularized estimators with some modifications for moderate sample sizes, and shed light on the influence of the alignment between the model parameter vector and the selected kernel.
\end{enumerate}

The remaining part of this paper is organized as follows. In Section \ref{sec: Empirical Bayes estimation}, we introduce preliminaries on the finite impulse response estimation problem and define the XMSE of an EB estimator. In Section \ref{sec:XMSE of empirical Bayes estimator}, we derive an explicit expression for the XMSE of an EB estimator with a general hyper-parameter estimator. In Section \ref{sec:accuracy of XMSE}, we introduce certain approximations to improve the accuracy of the XMSE expressions for regularized estimators for moderate sample sizes. In Section \ref{sec: numerical simulation}, we construct systems with positive XMSE of the EB-based regularized estimator for a finite sample size setting and show the average performance of different regularized estimators. In Section \ref{sec:conclusion}, we conclude this paper. The proofs of theorems and corollaries are included in Appendix A.

\textbf{Notation.} The set of all real-valued matrices of dimension $m_{1}\times m_{2}$ is denoted ${\R}^{m_{1}\times m_{2}}$. The $m$-by-$m$ dimensional identity matrix is denoted $\mathbf{I}_{m}$.
The first-order and second-order partial derivatives of $f(\bm{x}):\R^{m}\to \R$ with respect to $\bm{x}\in\R^{m}$ are denoted ${\partial f(\bm{x})}/{\partial \bm{x}}\in\R^{m}$ and ${\partial^2 f(\bm{x})}/{\partial \bm{x}\partial \bm{x}^{\top}}\in\R^{m\times m}$, respectively. The $(k,l)$th entry, the $k$th row and the $k$th column of a matrix $\mathbf{A}$ are $[\mathbf{A}]_{k,l}$, $[\mathbf{A}]_{k,:}$ and $[\mathbf{A}]_{:.k}$, respectively. The transpose and inverse of a matrix $\mathbf{A}$ are denoted $\mathbf{A}^{\top}$ and $\mathbf{A}^{-1}$, respectively. For a random vector $\bm{a}$, its expectation and covariance matrix are denoted $\E(\bm{a})$ and $\Var(\bm{a})\triangleq\E[(\bm{a}-\E(\bm{a}))(\bm{a}-\E(\bm{a}))^{\top}]$. The cross covariance matrix of random vectors $\bm{a}$ and $\bm{b}$ is denoted $\COV(\bm{a},\bm{b})\triangleq\E[(\bm{a}-\E(\bm{a}))(\bm{b}-\E(\bm{b}))^{\top}]$. A Gaussian distributed random variable $\bm\xi\in\R^{m}$ with mean $\bm{a}\in\R^{m}$ and covariance matrix $\mathbf{A}\in\R^{m\times m}$ is denoted $\mathcal{N}(\bm{a},\mathbf{A})$. The positive definite matrix $\mathbf{A}$ is denoted $\mathbf{A}\succ 0$. The rank of a matrix is denoted $\rank(\cdot)$. The Euclidean norm of a vector and the Frobenius norm of a matrix is denoted $\|\cdot\|_{2}$ and $\|\cdot\|_{F}$, respectively. The trace and determinant of a square matrix are denoted $\Tr(\cdot)$ and $\det(\cdot)$, respectively. A diagonal matrix $\mathbf{A}\in\R^{m\times m}$ with diagonal entries being one column vector $\bm{a}\in\R^{m}$ is denoted $\mathbf{A}=\diag\{\bm{a}\}$. Eigenvalues of a square matrix $\mathbf{A}\in\R^{m\times m}$ are denoted $\lambda_{1}(\mathbf{A})\geq \cdots\geq \lambda_{m}(\mathbf{A})$. The condition number of $\mathbf{A}$ is denoted $\cond(\mathbf{A})=\lambda_{1}(\mathbf{A})/\lambda_{m}(\mathbf{A})$.

\section{Empirical Bayes Estimation of Finite Impulse Response Parameters}\label{sec: Empirical Bayes estimation}

We first introduce the finite impulse response estimation problem and then define the corresponding EB estimators.

Consider an FIR model with scalar-valued input and output,
\begin{align}\label{eq:FIR model}
y(t)=\sum_{k=1}^{n}g_{k}u(t-k)+e(t),\quad t=1,\cdots,N,
\end{align}
where $y(t)\in\R$, $u(t)\in\R$ and $e(t)\in\R$ are the measurement output, input and measurement noise at time instant $t$, respectively, and $\{g_{k}\}_{k=1}^{n}$ are the impulse response parameters to be estimated, $n$ is the FIR model order, and $N$ is the sample size. For {notational} convenience, we reformulate \eqref{eq:FIR model} as
\begin{align}\label{eq:linear model}
\bm{Y}=\mathbf{\Phi}\bm\theta+\bm{E},
\end{align}
where $\bm{Y}=[y(1),\cdots,y(N)]^{\top}$, $\bm{E}=[e(1),\cdots,e(N)]^{\top}$, $\mathbf{\Phi}\in\R^{N\times n}$ is a lower triangular Toeplitz matrix {with $[\mathbf{\Phi}]_{:,1}=[u(0),\cdots,u(N-1)]^{\top}$} and known as the regression matrix, and $\bm\theta\in\R^{n}$ is the unknown parameter vector. We also need {the following} assumption.


\begin{assumption}\label{asp:input and noise}
\begin{enumerate}
\item[]
\item The FIR model order $n$ is fixed, and larger or equal to the true order.
\item The inputs $\{u(t)\}_{t=1-n}^{N-1}$ are deterministic and known with $u(t)=0$ for $t\leq 0$. Moreover, the regression matrix $\mathbf{\Phi}$ has full column rank for $N\geq n$, i.e., $\rank(\mathbf{\Phi})=n$.
\item The measurement noise is independent and identically Gaussian distributed with zero mean and known variance\footnote{In practice, the measurement noise variance can also be estimated from observations, see, e.g., Remark \ref{rmk:complete MSS of noise variance estimator} and \cite[Remark 5]{PDCDL14}.} $\sigma^2>0$, i.e., $\bm{E}\sim\mathcal{N}(\bm{0},\sigma^2\mathbf{I}_{N})$. 
\end{enumerate}
\end{assumption}

The performance of {an} estimator $\hat{\bm\theta}\in\R^{n}$ of the true model parameter $\bm\theta_{0}\in\R^{n}$ is measured by
\begin{align}\label{eq:def of MSE}
\MSE(\hat{\bm\theta})=\E[\|\hat{\bm\theta}-\bm\theta_{0}\|_{2}^2],
\end{align}
where the expectation is with respect to the measurement noise $\bm{E}$. 
{The objective of the finite impulse response estimation problem is to estimate the unknown $\bm\theta$ based on observed data $\{u(t),y(t)\}_{t=1}^{N}$}.

\subsection{Finite impulse response empirical Bayes estimators}\label{subsec: FIR parameter estimators}

For a squared loss function, the generalized Bayes estimator \cite{Lehmann:06} of $\bm\theta$ in \eqref{eq:linear model} is defined as
\begin{align}\label{eq:def of Bayes estimator}
\hat{\bm\theta}^{\Bayes}=&\argmin_{\hat{\bm\theta}\in\R^{n}}\int \|\hat{\bm\theta}-\bm\theta\|_{2}^2 \frac{p(\bm{Y}|\bm{\theta})\pi(\bm\theta)}{\int p(\bm{Y}|\bm{\theta})\pi(\bm\theta)d\bm\theta}d\bm\theta\nonumber\\
=&\frac{\int \bm\theta p(\bm{Y}|\bm{\theta})\pi(\bm\theta)d\bm\theta}{\int p(\bm{Y}|\bm{\theta})\pi(\bm\theta)d\bm\theta},
\end{align}
where $p(\bm{Y}|\bm\theta)$ is the probability density function (pdf) of $\bm{Y}|\bm\theta\sim\mathcal{N}(\mathbf{\Phi}\bm\theta,\sigma^2\mathbf{I}_{N})$, and $\pi(\bm\theta)$ is nonnegative and known as the weighting (prior) function. Note that  $\pi(\bm\theta)$ is allowed to be improper, i.e., $\int\pi(\bm\theta)d\bm\theta=\infty$. In particular, if $\pi(\bm\theta)=1$, the corresponding generalized Bayes estimator is equivalent to the ML estimator, given by
\begin{align}\label{eq:ML estimate}
\hat{\bm\theta}^{\ML}
=(\mathbf{\Phi}^{\top}\mathbf{\Phi})^{-1}\mathbf{\Phi}^{\top}\bm{Y}.
\end{align}


If the weighting function $\pi(\bm\theta)$ is selected from a family of weighting functions indexed by a hyper-parameter vector $\bm\eta\in\mathcal{D}_{\bm\eta}\subset\R^{p}$, where $\mathcal{D}_{\bm\eta}$ is the feasible set of $\bm\eta$, it is denoted as $\pi(\bm\theta|\bm\eta)$. The formal Bayes estimator \eqref{eq:def of Bayes estimator} is then generalized to the EB estimator by replacing $\pi(\bm\theta)$ with $\pi(\bm\theta|\hat{\bm\eta})$, where $\hat{\bm\eta}$ is estimated from observations $\{u(t),y(t)\}_{t=1}^{N}$.
More specifically, an EB estimator of $\bm\theta$ is defined as
\begin{align}
\hat{\bm\theta}^{\EB}(\hat{\bm\eta})
\label{eq: empirical Bayes estimator}
=&\frac{\int \bm\theta p(\bm{Y}|\bm{\theta})\pi(\bm\theta|\hat{\bm\eta})d\bm\theta}{\int p(\bm{Y}|\bm{\theta})\pi(\bm\theta|\hat{\bm\eta})d\bm\theta}.
\end{align}
The form of $\hat{\bm\theta}^{\EB}(\hat{\bm\eta})$ in \eqref{eq: empirical Bayes estimator} depends on the weighting function $\pi(\bm\theta|\bm\eta)$ and the hyper-parameter estimator $\hat{\bm\eta}$.

In particular, we consider the Gaussian weighting function, i.e., $\pi(\bm\theta|{\bm\eta})$ is in the form of the pdf of $\bm\theta|{\bm\eta}\sim\mathcal{N}(\bm{0},\mathbf{P}({\bm\eta}))$, where $\mathbf{P}(\bm\eta)\in\R^{n}$ is assumed to be positive definite. The corresponding EB estimator, referred to as the kernel-based regularized estimator, takes the form of
\begin{align}\label{eq:RLS estimate}
\hat{\bm\theta}^{\TR}({\bm\eta})
=&[\mathbf{\Phi}^{\top}\mathbf{\Phi}+\sigma^2\mathbf{P}({\bm\eta})^{-1}]^{-1}\mathbf{\Phi}^{\top}\bm{Y},
\end{align}
where the $(k,l)$th entry of the kernel matrix $\mathbf{P}(\bm\eta)$ is a specific kernel $[\mathbf{P}(\bm\eta)]_{k,l}=\kappa(k,l;\bm\eta)$. For the design of regularized estimators, we {consider} the following two problems.

The first problem is the selection of the kernel $\kappa(k,l;\bm\eta)$, which should leverage the given information of the true system to be estimated\cite{Chen18}. {A common option} is the stable spline (SS) kernel \cite{PN10a}, 
\begin{subequations}\label{eq:SS kernel}
\noeqref{eq:expression of SS kernel, eq:range of SS kernel}
\begin{align}\label{eq:expression of SS kernel}
\kappa_{SS}(k,l;\bm\eta)=&c\left(\frac{\gamma^{k+l+\max(k,l)}}{2}-\frac{\gamma^{3\max(k,l)}}{6}\right),\\
\label{eq:range of SS kernel}
\bm\eta=&[c,\gamma]^{\top},\ \mathcal{D}_{\bm\eta}=\{c\geq 0,\ 0\leq\gamma<1\}.\quad
\end{align}
\end{subequations}
It is parameterized by the two hyper-parameters $c$ and $\gamma$, which describe the scale and decay rate of the impulse response to be estimated, respectively. 

The second problem is how to estimate the hyper-parameter $\bm\eta$ from observed data $\{u(t),y(t)\}_{t=1}^{N}$. While our results apply to a general family of hyper-parameter estimators, we will provide explicit expressions for the following methods (see, e.g., \cite{PDCDL14, PCCDL22}) and their scaled generalizations. It is straightforward to derive explicit results for other estimators using Theorem \ref{thm:XMSE of empirical Bayes estimator} below.

\begin{itemize}
\item {\it Scaled} EB: The scaled EB hyper-parameter estimator is given by 
\begin{subequations}\label{eq: scaled EB}\noeqref{eq:EB alpha hyper-parameter estimator}\noeqref{eq:EB alpha cost function}
\begin{align}\label{eq:EB alpha hyper-parameter estimator}
\hat{\bm\eta}_{\EB,\alpha}=&\argmin_{\bm\eta\in\mathcal{D}_{\bm\eta}}\mathscr{F}_{\EB,\alpha}(\bm\eta),\\
\label{eq:EB alpha cost function}
\mathscr{F}_{\EB,\alpha}(\bm\eta)=&\bm{Y}^{\top}\mathbf{Q}(\bm\eta)^{-1}\bm{Y}+\alpha\log\det(\mathbf{Q}(\bm\eta)),
\end{align}
\end{subequations}
where $\mathbf{Q}(\bm\eta)=\mathbf{\Phi}\mathbf{P}(\bm\eta)\mathbf{\Phi}^{\top}+\sigma^2\mathbf{I}_{N}$. For convenience, $\hat{\bm\eta}_{\EB,\alpha}$ and $\mathscr{F}_{\EB,\alpha}(\bm\eta)$ with $\alpha=1$ will be denoted $\hat{\bm\eta}_{\EB}$ and $\mathscr{F}_{\EB}(\bm\eta)$, respectively. Notice that $\hat{\bm\eta}_{\EB}$ corresponds to the standard EB hyper-parameter estimator, which takes the pdf of $\bm\theta|\bm\eta \sim \mathcal{N}(0,\mathbf{P}(\bm\eta))$ as the weighting function $\pi(\bm\theta|\bm\eta)$ and then minimizes $\int p(\bm{Y}|\bm{\theta})\pi(\bm\theta|\bm\eta)d\bm\theta$ with respect to $\bm\eta$. The inclusion of the parameter $\alpha>0$, which is novel, allows us to strike another balance between the data fit and the model complexity terms in $\mathscr{F}_{\EB}(\bm\eta)$.

\item {\it Scaled Stein's unbiased risk estimation in relation to output prediction} ($\SURE_{y}$): We introduce the scaled $\SURE_{y}$ hyper-parameter estimator to be defined as
\begin{subequations}\label{eq: scaled SUREy}\noeqref{eq:SUREy alpha hyper-parameter estimator}\noeqref{eq: SUREy alpha cost function}
\begin{align}\label{eq:SUREy alpha hyper-parameter estimator}
&\hat{\bm\eta}_{\Sy,\alpha}=\argmin_{\bm\eta\in\mathcal{D}_{\bm\eta}}\mathscr{F}_{\Sy,\alpha}(\bm\eta),\\
\label{eq: SUREy alpha cost function}
&\mathscr{F}_{\Sy,\alpha}(\bm\eta)=\|\bm{Y}-\mathbf{H}(\bm\eta)\bm{Y}\|_{2}^2+2\alpha\sigma^2\Tr[\mathbf{H}(\bm\eta)],\quad 
\end{align}
\end{subequations}
where $\mathbf{H}(\bm\eta)=\mathbf{\Phi}\mathbf{P}(\bm\eta)\mathbf{\Phi}^{\top}\mathbf{Q}(\bm\eta)^{-1}$. We let $\hat{\bm\eta}_{\Sy}$ and $\mathscr{F}_{\Sy}(\bm\eta)$ denote $\hat{\bm\eta}_{\Sy,\alpha}$ and $\mathscr{F}_{\Sy,\alpha}(\bm\eta)$ with $\alpha=1$, respectively. Actually, $\mathscr{F}_{\Sy}(\bm\eta)$ is the Stein's unbiased risk estimator of the MSE in relation to output prediction, i.e., $\E[\mathscr{F}_{\Sy}(\bm\eta)]={\MSE}_{\ty}(\hat{\bm\theta}^{\TR}(\bm\eta))$, where 
\begin{align}
\label{eq:def of MSEy}\noeqref{eq:def of MSEy}
{\MSE}_{\ty}(\hat{\bm\theta}^{\TR}(\bm\eta))=&\E[\|\mathbf{\Phi}\bm\theta_{0}+\widetilde{\bm{E}}-\mathbf{\Phi}\hat{\bm\theta}^{\TR}(\bm\eta)\|_{2}^2],
\end{align}
with $\widetilde{\bm{E}}\sim\mathcal{N}(\bm{0},\sigma^2 \mathbf{I}_{N})$ independent of $\bm{E}$ in \eqref{eq:linear model}. Thus, $\hat{\bm\eta}_{\Sy}$ corresponds to the standard $\SURE_{y}$ hyper-parameter estimator. The parameter $\alpha>0$ allows us to alter the balance between the variance and bias in the estimator.
\item {\it Scaled} GCV: The scaled GCV hyper-parameter estimator can be expressed as
\begin{subequations}
\label{eq: scaled GCV}\noeqref{eq:GCV alpha hyper-parameter estimator}\noeqref{eq:GCV alpha cost function}
\begin{align}
\label{eq:GCV alpha hyper-parameter estimator}
&\hat{\bm\eta}_{\GCV,\alpha}=\argmin_{\bm\eta\in\mathcal{D}_{\bm\eta}}\mathscr{F}_{\GCV,\alpha}(\bm\eta)\\
\label{eq:GCV alpha cost function}
&\mathscr{F}_{\GCV,\alpha}(\bm\eta)=\|\bm{Y}-\mathbf{H}(\bm\eta)\bm{Y}\|_{2}^{2}\nonumber\\
&+\alpha\|\bm{Y}-\mathbf{H}(\bm\eta)\bm{Y}\|_{2}^{2}\left\{\frac{1}{\left[1-\Tr(\mathbf{H}(\bm\eta))/N\right]^2}-1\right\}.\quad\ 
\end{align}
\end{subequations}
For $\alpha=1$, $\hat{\bm\eta}_{\GCV,\alpha}$ and $\mathscr{F}_{\GCV,\alpha}(\bm\eta)$ will be denoted $\hat{\bm\eta}_{\GCV}$ and $\mathscr{F}_{\GCV}(\bm\eta)$, respectively. Note that $\hat{\bm\eta}_{\GCV}$ corresponds to the standard GCV hyper-parameter estimator and is a simplification of the predicted residual error sum of squares (PRESS) hyper-parameter estimator. Similarly, the parameter $\alpha>0$ is to modify the bias-variance trade-off of $\hat{\bm\eta}_{\GCV}$.
\end{itemize}

\subsection{Motivating example}\label{subsec:motivating example}

As discussed in \cite{COL12a,PDCDL14}, the regularized estimator $\hat{\bm\theta}^{\TR}(\hat{\bm\eta})$ in \eqref{eq:RLS estimate} can reduce the possibly large variance of the ML estimator $\hat{\bm\theta}^{\ML}$ in \eqref{eq:ML estimate} by introducing some bias, especially for small sample sizes $N$ or low signal-to-noise ratio ({SNR}).

\begin{figure}[!htbp]
\centering
\includegraphics[width=0.9\linewidth]{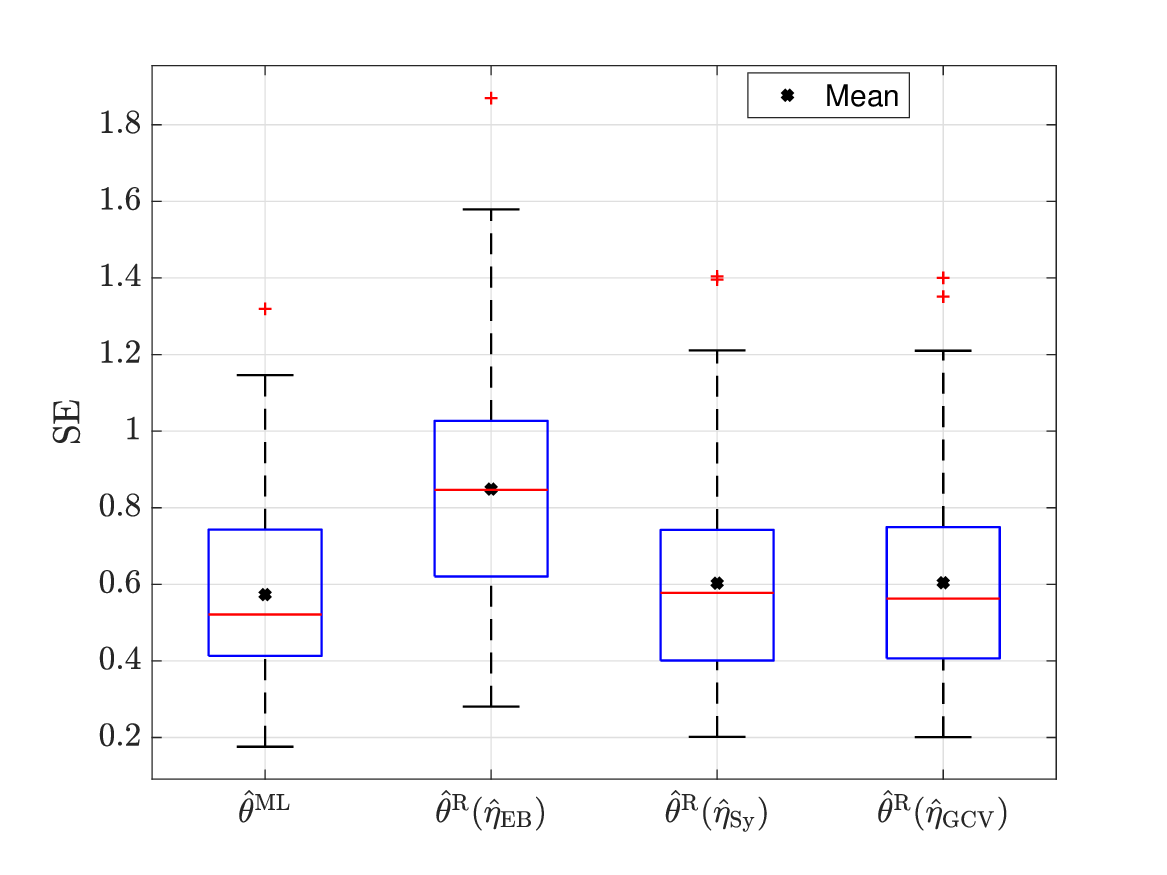}
\caption{Boxplots of the SE of $\hat{\bm\theta}^{\ML}$, $\hat{\bm\theta}^{\TR}(\hat{\eta}_{\EB})$, $\hat{\bm\theta}^{\TR}(\hat{\eta}_{\Sy})$ and $\hat{\bm\theta}^{\TR}(\hat{\eta}_{\GCV})$ for one specific system. Their sample MSE are $5.73\times 10^{-1}$, $8.49\times 10^{-1}$, $6.03\times 10^{-1}$ and $6.04\times 10^{-1}$, respectively.}
\label{fig:motivating example}
\end{figure}

However, let us now consider a system generated in Matlab by \lstinline[style=Matlab-editor]{rss(30,1,1)}{and the method} in \cite[Section 2]{COL12a} such that the largest pole modulus is smaller than $0.95$. Its impulse response is truncated at $n=20$ to obtain $\tilde{\bm\theta}_{0}\in\R^{20}$, which {contains} the first $20$ impulse response coefficients. We generate inputs $\{u(t)\}_{t=1}^{N}$ as independent realizations of Gaussian noise with zero mean and unit variance, and use sample size $N=50$ and measurement noise variance $\sigma^2=1$. Furthermore, we scale $\bm\theta_{0}=m\tilde{\bm\theta}_{0}$ and select $m\in\R$ such that the sample SNR, i.e., the ratio between the sample variance of the noiseless outputs $\mathbf{\Phi}\bm\theta_{0}$ and $\sigma^2=1$, is $10$. For fixed $\bm\theta_{0}$ and $\{u(t)\}_{t=1}^{N}$, we run $100$ Monte Carlo (MC) simulations to obtain $100$ output sequences $\{y(t)\}_{t=1}^{N}$. For each MC simulation, we calculate the squared error (SE) of $\hat{\bm\theta}$,
\begin{align}\label{eq:squared error of theta_hat}
\text{SE}(\hat{\bm\theta})=\|\hat{\bm\theta}-\bm\theta_{0}\|_{2}^2.
\end{align}
The sample mean of $\text{SE}(\hat{\bm\theta})$ over the MC simulations is referred to as the ``sample $\MSE(\hat{\bm\theta})$''. In Fig. \ref{fig:motivating example}, we show boxplots for the SE of $\hat{\bm\theta}^{\ML}$, $\hat{\bm\theta}^{\TR}(\hat{\eta}_{\EB})$, $\hat{\bm\theta}^{\TR}(\hat{\eta}_{\Sy})$ and $\hat{\bm\theta}^{\TR}(\hat{\eta}_{\GCV})$. Here, we select the SS kernel \eqref{eq:SS kernel} with $\gamma=0.95$ for all regularized estimators. It can be observed that $\hat{\bm\theta}^{\ML}$ outperforms all three regularized estimators. This observation seems ``contradictory'' to the previous understanding of regularized estimators. It motivates us to analyze the MSE of the regularized estimator in \eqref{eq:RLS estimate}, or more generally, the MSE of the EB estimator in \eqref{eq: empirical Bayes estimator}.

\begin{remark}\label{rmk:when regularization cannot improve}
Note that the performance comparison between the ML and kernel-based regularized estimators has been partially discussed in \cite{MLC24}. It concludes that when the regression matrix $\mathbf{\Phi}$ and the kernel matrix $\mathbf{P}(\hat{\bm\eta})$ are well-conditioned, the probability that $\mathbf{P}(\hat{\bm\eta})=+\infty$, leading to $\hat{\bm\theta}^{\TR}(\hat{\bm\eta})|_{\mathbf{P}(\hat{\bm\eta})=+\infty}=\hat{\bm\theta}^{\ML}$, is significantly larger than zero. This explains when the performance of regularized estimators is the same as that of the ML estimator. However, the hyper-parameter estimates in the motivating example are not excessively large, and $\mathbf{P}(\hat{\bm\eta})$ does not approach $+\infty$. Therefore, the motivating example cannot be explained by the analysis in \cite{MLC24}.
\end{remark}

\subsection{Analyzing MSE of an empirical Bayes estimator}

To understand the observation in the motivating example, we need an analytic criterion to analyze the performance of an EB estimator in \eqref{eq: empirical Bayes estimator}. It is well-known that $\MSE(\hat{\bm\theta}^{\ML})=\sigma^2\Tr[(\mathbf{\Phi}^{\top}\mathbf{\Phi})^{-1}]$. However, the MSE of $\hat{\bm\theta}^{\EB}(\hat{\bm\eta})$ is generally the expectation of a complicated nonlinear function and has no explicit form for finite sample sizes. For this, we shift our focus to large sample sizes and add the following assumption regarding the external excitation.

\begin{assumption}\label{asp:limit of PPN}
The regression matrix $\mathbf{\Phi}$ in \eqref{eq:FIR model} satisfies
\begin{align}\label{eq:limit of PPN}
\lim_{N\to\infty}\frac{\mathbf{\Phi}^{\top}\mathbf{\Phi}}{N}=\mathbf{\Sigma}\succ 0.
\end{align}
\end{assumption}

We start by deriving the limits of $\MSE(\hat{\bm\theta}^{\ML})$, $\MSE(\hat{\bm\theta}^{\TR}({\bm\eta}))$ and their difference, where $\hat{\bm\theta}^{\TR}(\bm\eta)$ with fixed $\bm\eta$ can be seen as a specific EB estimator as discussed before \eqref{eq:RLS estimate}.

\begin{lemma}\label{lemma:limits of MSEs of ML and RLS estimators}
Under Assumptions \ref{asp:input and noise}-\ref{asp:limit of PPN}, for fixed $\bm\eta$, we have
\begin{subequations}
\begin{gather}
\label{eq:1st order limit of MSE of ML estimator}
\lim_{N\to\infty}N\MSE(\hat{\bm\theta}^{\ML})=\sigma^2\mathbf{\Sigma}^{-1},\\
\label{eq:1st order limit of MSE of RLS estimator}
\lim_{N\to\infty}N\MSE(\hat{\bm\theta}^{\TR}(\bm\eta))=\sigma^2\mathbf{\Sigma}^{-1},\\
\lim_{N\to\infty}N^2[\MSE(\hat{\bm\theta}^{\TR}(\bm\eta))-\MSE(\hat{\bm\theta}^{\ML})]=\nonumber\\
\label{eq:2nd order limit of MSE difference}
(\sigma^2)^2\|\mathbf{\Sigma}^{-1}\mathbf{P}(\bm\eta)^{-1}\bm\theta_{0}\|_{2}^2-2(\sigma^2)^2\Tr[\mathbf{\Sigma}^{-2}\mathbf{P}(\bm\eta)^{-1}].\ 
\end{gather}
\end{subequations}
\end{lemma}

As given in \eqref{eq:1st order limit of MSE of ML estimator}-\eqref{eq:1st order limit of MSE of RLS estimator}, the first-order limits of $\MSE(\hat{\bm\theta}^{\ML})$ and $\MSE(\hat{\bm\theta}^{\TR}({\bm\eta}))$ are the same. Thus, the first-order limit is inadequate to expose the difference between $\hat{\bm\theta}^{\ML}$ and $\hat{\bm\theta}^{\TR}({\bm\eta})$. We then discard the common term $\sigma^2\mathbf{\Sigma}^{-1}$ and consider the second-order limit of the difference between $\MSE(\hat{\bm\theta}^{\ML})$ and $\MSE(\hat{\bm\theta}^{\TR}({\bm\eta}))$ in \eqref{eq:2nd order limit of MSE difference}. This motivates us to consider a general EB estimator $\hat{\bm\theta}^{\EB}(\hat{\bm\eta})$
%
%
%
and define its excess MSE (XMSE) as
\begin{align}\label{eq:def of XMSE}
&\XMSE(\hat{\bm\theta}^{\EB}(\hat{\bm\eta}))\nonumber\\
=&\lim_{N\to\infty}N^2[\MSE(\hat{\bm\theta}^{\EB}(\hat{\bm\eta}))-\MSE(\hat{\bm\theta}^{\ML})],
\end{align}
for which more assumptions on the hyper-parameter estimator $\hat{\bm\eta}$ are needed and will be discussed shortly later. The XMSE of $\hat{\bm\theta}^{\EB}(\hat{\bm\eta})$ in \eqref{eq:def of XMSE} quantifies the performance difference between the EB and ML estimators for large sample sizes. If $\hat{\bm\theta}^{\ML}$ outperforms $\hat{\bm\theta}^{\EB}(\hat{\bm\eta})$, $\XMSE(\hat{\bm\theta}^{\EB}(\hat{\bm\eta}))$ is positive; otherwise, it is negative.


In what follows, we will first derive an explicit expression for the XMSE of an EB estimator in \eqref{eq: empirical Bayes estimator} using an arbitrary data-dependent hyper-parameter $\hat{\bm\eta}$. We then consider the generalized Bayes estimator, and the kernel-based regularized estimators with specific hyper-parameter estimators, including $\hat{\bm\eta}_{\EB,\alpha}$, $\hat{\bm\eta}_{\Sy,\alpha}$ and $\hat{\bm\eta}_{\GCV,\alpha}$, detailed in Section \ref{subsec: FIR parameter estimators}, and derive expressions for their XMSE. This will allow us to provide a theoretical analysis of the performance of EB estimators. We will also discuss the accuracy of the XMSE for capturing the difference in performance between ML and regularized estimators for moderate sample sizes.

\section{Excess MSE of an empirical Bayes estimator}\label{sec:XMSE of empirical Bayes estimator}

We {start} by observing that it is only of interest to study data-dependent hyper-parameters that are functions of the ML estimator $\hat{\bm\theta}^{\ML}$, denoted as $\hat{\bm\eta}(\hat{\bm\theta}^{\ML})$.
To understand why this is the case, we notice that the Rao-Blackwell theorem \cite[Theorem 1.7.8]{Lehmann:06} gives that if there exists a minimal sufficient statistic (MSS) $\bm{S}$, then $\MSE(\E(\hat{\bm\theta}|\bm{S}))\leq \MSE(\hat{\bm\theta})$, with equality if and only if $\hat{\bm\theta}=\E(\hat{\bm\theta}|\bm{S})$, for any estimator $\hat{\bm\theta}$. The implication is that it is only of interest to consider estimators based on an MSS, provided that such a quantity exists. Now for the Gaussian setting \eqref{eq:linear model}, $\hat{\bm\theta}^{\ML}$ is such a statistic \cite{Lehmann:06}, and hence only estimators of the type $\hat{\bm\theta}(\hat{\bm\theta}^{\ML})$ need to be considered. Since $\hat{\bm\theta}^{\ML}$ is an MSS, $\hat{\bm\theta}^{\EB}(\bm\eta)$ in \eqref{eq: empirical Bayes estimator} with fixed $\bm\eta$ is a function of $\hat{\bm\theta}^{\ML}$. Clearly, the EB estimator $\hat{\bm\theta}^{\EB}(\hat{\bm\eta}(\hat{\bm\theta}^{\ML}))$ is also a function of $\hat{\bm\theta}^{\ML}$.

\begin{remark}\label{rmk:functions of ML}
Note that $\hat{\bm\theta}^{\TR}(\bm\eta)$ in \eqref{eq:RLS estimate}, $\hat{\bm\eta}_{\EB,\alpha}$ in \eqref{eq:EB alpha hyper-parameter estimator}, $\hat{\bm\eta}_{\Sy,\alpha}$ in \eqref{eq:SUREy alpha hyper-parameter estimator} and $\hat{\bm\eta}_{\GCV,\alpha}$ in \eqref{eq:GCV alpha hyper-parameter estimator}, can all be reformulated as functions of $\hat{\bm\theta}^{\ML}$:
\begin{align}
\label{eq:reformulation of theta_RLS}
\hat{\bm\theta}^{\TR}(\bm\eta)=&\hat{\bm\theta}^{\ML}-\sigma^2(\mathbf{\Phi}^{\top}\mathbf{\Phi})^{-1}\mathbf{S}(\bm\eta)^{-1}\hat{\bm\theta}^{\ML},\\
\mathscr{F}_{\EB,\alpha}(\bm\eta)=&\mathscr{F}_{\tb,\alpha}(\hat{\bm\theta}^{\ML},\bm\eta)+\text{terms}\ \text{independent}\ \text{of}\ \bm\eta,\nonumber\\
\label{eq:reformulation of F_Sy}
\mathscr{F}_{\Sy,\alpha}(\bm\eta)=&\mathscr{F}_{\ty,\alpha}(\hat{\bm\theta}^{\ML},\bm\eta)+\text{terms}\ \text{independent}\ \text{of}\ \bm\eta,\\
\label{eq:reformulation of F_GCV}
\mathscr{F}_{\GCV,\alpha}(\bm\eta)=&\mathscr{F}_{\ty,\alpha}(\hat{\bm\theta}^{\ML},\bm\eta)+o(\|(\mathbf{\Phi}^{\top}\mathbf{\Phi})^{-1}\|_{F}),
\end{align}
where
\begin{align}
\label{eq:def of S}
&\mathbf{S}(\bm\eta)=\mathbf{P}(\bm\eta)+\sigma^2(\mathbf{\Phi}^{\top}\mathbf{\Phi})^{-1},\\
&\mathscr{F}_{\tb,\alpha}(\hat{\bm\theta}^{\ML},\bm\eta)=(\hat{\bm\theta}^{\ML})^{\top}\mathbf{S}(\bm\eta)^{-1}\hat{\bm\theta}^{\ML}+\alpha\log\det(\mathbf{S}(\bm\eta)),\nonumber\\
\label{eq:F_ty_alpha}
&\mathscr{F}_{\ty,\alpha}(\hat{\bm\theta}^{\ML},\bm\eta)=(\sigma^2)^2(\hat{\bm\theta}^{\ML})^{\top}\mathbf{S}(\bm\eta)^{-1}(\mathbf{\Phi}^{\top}\mathbf{\Phi})^{-1}\nonumber\\
&\times\mathbf{S}(\bm\eta)^{-1}\hat{\bm\theta}^{\ML}
-2\alpha(\sigma^2)^2\Tr[(\mathbf{\Phi}^{\top}\mathbf{\Phi})^{-1}\mathbf{S}(\bm\eta)^{-1}].
\end{align}
\end{remark}

\begin{remark}\label{rmk:complete MSS of noise variance estimator}
When the measurement noise variance $\sigma^2$ in Assumption \ref{asp:input and noise} is unknown, we typically estimate it by $\widehat{\sigma^2}={\|\bm{Y}-\mathbf{\Phi}\hat{\bm\theta}^{\ML}\|_{2}^2}/{(N-n)}$. In this case, $(\hat{\bm\theta}^{\ML},\widehat{\sigma^2})$ is a complete MSS, and thus it is necessary to consider hyper-parameter estimators $\hat{\bm\eta}$ that are functions of both $\hat{\bm\theta}^{\ML}$ and $\widehat{\sigma^2}$, i.e., $\hat{\bm\eta}(\hat{\bm\theta}^{\ML},\widehat{\sigma^2})$.
\end{remark}

\begin{remark}\label{rmk:complete MSS of arx}
For the autoregressive with external input (ARX) model, the complete MSS is $(\hat{\bm\theta}^{\ML},(\mathbf{\Phi}^{\top}\mathbf{\Phi})^{-1})$, and thus one has to consider hyper-parameter estimators of the form $\hat{\bm\eta}(\hat{\bm\theta}^{\ML},(\mathbf{\Phi}^{\top}\mathbf{\Phi})^{-1})$.
\end{remark}

Before {proceeding to} the XMSE of $\hat{\bm\theta}^{\EB}(\hat{\bm\eta}(\hat{\bm\theta}^{\ML}))$, we first decompose the difference between $\MSE(\hat{\bm\theta}^{\EB}(\hat{\bm\eta}(\hat{\bm\theta}^{\ML})))$ and {$\MSE(\hat{\bm\theta}^{\ML})$} in the following theorem.

\begin{theorem}\label{thm:decomposition of MSE difference}
Under Assumption \ref{asp:input and noise}, we have
\begin{align*}
&\MSE(\hat{\bm\theta}^{\EB}(\hat{\bm\eta}(\hat{\bm\theta}^{\ML})))-\MSE(\hat{\bm\theta}^{\ML})\\
=&\|\mathbf{\Upsilon}_{\Bias}\|_{2}^2+\Tr[\mathbf{\Upsilon}_{\Var}]+\Tr[\mathbf{\Upsilon}_{\VarHPE}]+\Tr[\mathbf{\Upsilon}_{\HOT}],
\end{align*}
where 
\begin{align}
&\mathbf{\Upsilon}_{\Bias}=\E[\hat{\bm\theta}^{\EB}(\hat{\bm\eta}(\bm\theta_{0}))]-\bm\theta_{0},\\
&\hat{\bm\eta}(\bm\theta_{0})=\argmin_{\bm\eta\in\mathcal{D}_{\bm\eta}}\mathscr{F}(\bm\theta_{0},\bm\eta),\\
&\mathbf{\Upsilon}_{\Var}=\COV(\hat{\bm\theta}^{\EB}(\hat{\bm\eta}(\bm\theta_{0}))-\hat{\bm\theta}^{\ML},\hat{\bm\theta}^{\ML})\\
&\quad +\COV(\hat{\bm\theta}^{\ML},\hat{\bm\theta}^{\EB}(\hat{\bm\eta}(\bm\theta_{0}))-\hat{\bm\theta}^{\ML}),\\
&\mathbf{\Upsilon}_{\VarHPE}=\COV(\hat{\bm\theta}^{\EB}(\hat{\bm\eta}(\hat{\bm\theta}^{\ML}))-\hat{\bm\theta}^{\EB}(\hat{\bm\eta}(\bm\theta_{0})),\hat{\bm\theta}^{\ML})\\
&\quad +\COV(\hat{\bm\theta}^{\ML},\hat{\bm\theta}^{\EB}(\hat{\bm\eta}(\hat{\bm\theta}^{\ML}))-\hat{\bm\theta}^{\EB}(\hat{\bm\eta}(\bm\theta_{0}))),\\
&\mathbf{\Upsilon}_{\HOT}=\Var(\hat{\bm\theta}^{\EB}(\hat{\bm\eta}(\hat{\bm\theta}^{\ML}))-\hat{\bm\theta}^{\ML})\nonumber\\
&\quad +[\E(\hat{\bm\theta}^{\EB}(\hat{\bm\eta}(\hat{\bm\theta}^{\ML})))-\bm\theta_{0}][\E(\hat{\bm\theta}^{\EB}(\hat{\bm\eta}(\hat{\bm\theta}^{\ML})))-\bm\theta_{0}]^{\top}\nonumber\\
&\quad -[\E[\hat{\bm\theta}^{\EB}(\hat{\bm\eta}(\bm\theta_{0}))]-\bm\theta_{0}][\E[\hat{\bm\theta}^{\EB}(\hat{\bm\eta}(\bm\theta_{0}))]-\bm\theta_{0}]^{\top}.
\end{align}
\end{theorem}

As defined in \eqref{eq:def of XMSE}, the expression of $\XMSE(\hat{\bm\theta}^{\EB}(\hat{\bm\eta}(\hat{\bm\theta}^{\ML})))$ depends on the limits of $N\mathbf{\Upsilon}_{\Bias}$, $N^2\mathbf{\Upsilon}_{\Var}$, $N^2\mathbf{\Upsilon}_{\VarHPE}$ and $N^2\mathbf{\Upsilon}_{\HOT}$. To derive these limits, we will apply Taylor expansions of $\hat{\bm\theta}^{\EB}(\hat{\bm\eta}(\hat{\bm\theta}^{\ML}))$ and $\hat{\bm\theta}^{\EB}(\hat{\bm\eta}(\bm\theta_{0}))$ at $\hat{\bm\theta}^{\ML}=\bm\theta_{0}$. Then, we need the limits of $\hat{\bm\eta}(\hat{\bm\theta}^{\ML})$ and its derivative given $\hat{\bm\theta}^{\ML}=\bm\theta_{0}$, which are provided in the following assumption.

\begin{assumption}\label{asp:limits of hyper-parameter estimator and its derivative}
{The following limits exist,}
\begin{align}\label{eq:convergence of hyper-parameter estimator}
\lim_{N\to\infty}\hat{\bm\eta}(\hat{\bm\theta}^{\ML})|_{\hat{\bm\theta}^{\ML}=\bm\theta_{0}}=&\bm\eta_{\star}(\bm\theta_{0}),\\
\label{eq:limit of derivative of hateta}
\lim_{N\to\infty}\left.\frac{\partial \hat{\bm\eta}(\hat{\bm\theta}^{\ML})}{\partial \hat{\bm\theta}^{\ML}}\right|_{\hat{\bm\theta}^{\ML}=\bm\theta_{0}}=&\mathbf{D}_{\bm\eta,\star}^{'}(\bm\theta_{0}).
\end{align}
\end{assumption}


Based on Assumption \ref{asp:limits of hyper-parameter estimator and its derivative}, we establish the limits of $N\mathbf{\Upsilon}_{\Bias}$, $N^2\mathbf{\Upsilon}_{\Var}$, $N^2\mathbf{\Upsilon}_{\VarHPE}$ and $N^2\mathbf{\Upsilon}_{\HOT}$ in Theorem \ref{thm:decomposition of MSE difference}, and provide an expression for the XMSE of $\hat{\bm\theta}^{\EB}(\hat{\bm\eta}(\hat{\bm\theta}^{\ML}))$.

\begin{theorem}\label{thm:XMSE of empirical Bayes estimator}
Under Assumptions \ref{asp:input and noise}-\ref{asp:limits of hyper-parameter estimator and its derivative}, we have
\begin{gather*}
\lim_{N\to\infty}N\mathbf{\Upsilon}_{\Bias}=\XBias(\hat{\bm\theta}^{\EB}(\hat{\bm\eta}(\hat{\bm\theta}^{\ML}))),\\
\lim_{N\to\infty}N^2\mathbf{\Upsilon}_{\Var}=\XVar(\hat{\bm\theta}^{\EB}(\hat{\bm\eta}(\hat{\bm\theta}^{\ML}))),\\
\lim_{N\to\infty}N^2\mathbf{\Upsilon}_{\VarHPE}=\XVarHPE(\hat{\bm\theta}^{\EB}(\hat{\bm\eta}(\hat{\bm\theta}^{\ML})))
\end{gather*}
and $\lim_{N\to\infty}N^2\mathbf{\Upsilon}_{\HOT}=\mathbf{0}$, where
\begin{subequations}
\begin{align}\label{eq:def of XBias}
&\XBias(\hat{\bm\theta}^{\EB}(\hat{\bm\eta}(\hat{\bm\theta}^{\ML})))
=\bm{b}_{\bm\theta_{0},\star}(\bm\eta_{\star}(\bm\theta_{0})),\\
\label{eq:def of XVar}
&\XVar(\hat{\bm\theta}^{\EB}(\hat{\bm\eta}(\hat{\bm\theta}^{\ML})))=\nonumber\\
&\quad{\XVar}_{\#}(\hat{\bm\theta}^{\EB}(\hat{\bm\eta}(\hat{\bm\theta}^{\ML})))
+{\XVar}_{\#}(\hat{\bm\theta}^{\EB}(\hat{\bm\eta}(\hat{\bm\theta}^{\ML})))^{\top},\\
\label{eq:def of XVar part1}
&{\XVar}_{\#}(\hat{\bm\theta}^{\EB}(\hat{\bm\eta}(\hat{\bm\theta}^{\ML})))=\sigma^2\mathbf{b}_{\bm\theta_{0},\star}^{'}(\bm\eta_{\star}(\bm\theta_{0}))\mathbf{\Sigma}^{-1},\\
\label{eq:def of XVarHPE}
&\XVarHPE(\hat{\bm\theta}^{\EB}(\hat{\bm\eta}(\hat{\bm\theta}^{\ML})))=
{\XVarHPE}_{\#}(\hat{\bm\theta}^{\EB}(\hat{\bm\eta}(\hat{\bm\theta}^{\ML})))\nonumber\\
&\quad+{\XVarHPE}_{\#}(\hat{\bm\theta}^{\EB}(\hat{\bm\eta}(\hat{\bm\theta}^{\ML})))^{\top},\\
\label{eq:def of VarHPE part1}
&{\XVarHPE}_{\#}(\hat{\bm\theta}^{\EB}(\hat{\bm\eta}(\hat{\bm\theta}^{\ML})))=\nonumber\\
&\qquad\sigma^2\mathbf{b}_{\bm\eta,\star}^{'}(\bm\eta)|_{\bm\eta=\bm\eta_{\star}(\bm\theta_{0})}\mathbf{D}_{\bm\eta,\star}^{'}(\bm\theta_{0})\mathbf{\Sigma}^{-1},\\
\label{eq:diff of b at theta0}
&\bm{b}_{\bm\theta_{0},\star}(\bm\eta)=\lim_{N\to\infty}N[\hat{\bm\theta}^{\EB}(\bm\eta)|_{\hat{\bm\theta}^{\ML}=\bm\theta_{0}}-\bm\theta_{0}]\nonumber\\
&\qquad\quad\ =\sigma^2\mathbf{\Sigma}^{-1}\left.\frac{\partial \log(\pi(\bm\theta|\bm\eta))}{\partial \bm\theta}\right|_{\bm\theta=\bm\theta_{0}},\\
\label{eq:derivative of b wrt theta}
&\mathbf{b}_{\bm\theta_{0},\star}^{'}(\bm\eta)=\lim_{N\to\infty}N\left[\left.\frac{\partial \hat{\bm\theta}^{\EB}(\bm\eta)}{\partial \hat{\bm\theta}^{\ML}}\right|_{\hat{\bm\theta}^{\ML}=\bm\theta_{0}}-\mathbf{I}_{n}\right]\nonumber\\
&\qquad\quad\ \ =\sigma^2\mathbf{\Sigma}^{-1}\left.\frac{\partial^2 \log(\pi(\bm\theta|\bm\eta))}{\partial \bm\theta\partial \bm\theta^{\top}}\right|_{\bm\theta=\bm\theta_{0}},\\
\label{eq:derivative of b wrt eta}
&\mathbf{b}_{\bm\eta,\star}^{'}(\bm\eta)=\lim_{N\to\infty}N\left.\frac{\partial \hat{\bm\theta}^{\EB}(\bm\eta)}{\partial \bm\eta}\right|_{\hat{\bm\theta}^{\ML}=\bm\theta_{0}}\nonumber\\
&\qquad\quad\  =\sigma^2\mathbf{\Sigma}^{-1}\left.\frac{\partial^2 \log(\pi(\bm\theta|\bm\eta))}{\partial \bm\theta\partial \bm\eta^{\top}}\right|_{\bm\theta=\bm\theta_{0}}.
\end{align}
\end{subequations}
Moreover, the XMSE of $\hat{\bm\theta}^{\EB}(\hat{\bm\eta}(\hat{\bm\theta}^{\ML}))$ equals
\begin{align}\label{eq:expression of XMSE}
\XMSE(\hat{\bm\theta}^{\EB}(\hat{\bm\eta}(\hat{\bm\theta}^{\ML})))&=\|\XBias(\hat{\bm\theta}^{\EB}(\hat{\bm\eta}(\hat{\bm\theta}^{\ML})))\|_{2}^2\nonumber\\
+\Tr[&\XVar(\hat{\bm\theta}^{\EB}(\hat{\bm\eta}(\hat{\bm\theta}^{\ML})))]\nonumber\\
+\Tr[&\XVarHPE(\hat{\bm\theta}^{\EB}(\hat{\bm\eta}(\hat{\bm\theta}^{\ML})))].
\end{align}
\end{theorem}

\begin{remark}\label{rmk:improper weighting}
Even if the weighting function $\pi(\bm\theta|\hat{\bm\eta}(\hat{\bm\theta}^{\ML}))$ is improper, the definition of $\hat{\bm\theta}^{\EB}(\hat{\bm\eta}(\hat{\bm\theta}^{\ML}))$ in \eqref{eq: empirical Bayes estimator} contains an implicit condition that 
\begin{align}\label{eq:two proper terms}
\bm\theta\pi(\bm\theta|\hat{\bm\eta}(\hat{\bm\theta}^{\ML}))p(\bm{Y}|\bm\theta)\ \text{and}\ \pi(\bm\theta|\hat{\bm\eta}(\hat{\bm\theta}^{\ML}))p(\bm{Y}|\bm\theta)
\end{align}
are both well-defined and integrable. We can see from Section \ref{subsec: proof of XMSE of empirical Bayes estimator} that the derivation of the XMSE of $\hat{\bm\theta}^{\EB}(\hat{\bm\eta}(\hat{\bm\theta}^{\ML}))$ only involves the Taylor expansions of \eqref{eq:two proper terms} and their integrals. Hence, Theorem \ref{thm:XMSE of empirical Bayes estimator} remains true for improper $\pi(\bm\theta|\hat{\bm\eta}(\hat{\bm\theta}^{\ML}))$.
\end{remark}

The expression \eqref{eq:expression of XMSE} for the XMSE of $\hat{\bm\theta}^{\EB}(\hat{\bm\eta}(\hat{\bm\theta}^{\ML}))$ contains three components. Note that $\MSE(\hat{\bm\theta}^{\EB}(\hat{\bm\eta}(\hat{\bm\theta}^{\ML})))$ is non-negative, while $\XMSE(\hat{\bm\theta}^{\EB}(\hat{\bm\eta}(\hat{\bm\theta}^{\ML}))))$ is the asymptotic MSE difference between two estimators, which is not necessarily non-negative. The sign of $\XMSE(\hat{\bm\theta}^{\EB}(\hat{\bm\eta}(\hat{\bm\theta}^{\ML}))))$ is determined by its three components:
\begin{itemize}
\item[-] The first component $\|\XBias(\hat{\bm\theta}^{\EB}(\hat{\bm\eta}(\hat{\bm\theta}^{\ML})))\|_{2}^2$ is denoted the excess squared bias of $\hat{\bm\theta}^{\EB}(\hat{\bm\eta}(\hat{\bm\theta}^{\ML}))$, and clearly is non-negative. Note that \eqref{eq:def of XBias} is the limit of $N[\E(\hat{\bm\theta}^{\EB}(\hat{\bm\eta}(\bm\theta_{0})))-\bm\theta_{0}]$, and thus {the first component describes the squared asymptotic bias} of $\hat{\bm\theta}^{\EB}(\hat{\bm\eta}(\bm\theta_{0}))$.
\item[-] The second component $\Tr[\XVar(\hat{\bm\theta}^{\EB}(\hat{\bm\eta}(\hat{\bm\theta}^{\ML})))]$ is {denoted} the excess variance caused by $\hat{\bm\theta}^{\EB}(\hat{\bm\eta}(\bm\theta_{0}))$. This term thus captures the excess variance when the hyper-parameter is kept fixed to its limiting value. The sign of 
\begin{align*}
&\Tr[\XVar(\hat{\bm\theta}^{\EB}(\hat{\bm\eta}(\hat{\bm\theta}^{\ML})))]\\
=&2(\sigma^2)^2\Tr\left[\mathbf{\Sigma}^{-1}\left.\frac{\partial^2 \log(\pi(\bm\theta|\bm\eta_{\star}(\bm\theta_{0})))}{\partial \bm\theta\partial\bm\theta^{\top}}\right|_{\bm\theta=\bm\theta_{0}}\mathbf{\Sigma}^{-1}\right]
\end{align*}
depends on the form of $\pi(\bm\theta|\bm\eta)$. More specifically, it depends on the curvature of $\log(\pi(\bm\theta|\bm\eta))$ at the convergence point of $\hat{\bm\eta}(\hat{\bm\theta}^{\ML})$.

\item[-] The third component $\Tr[\XVarHPE(\hat{\bm\theta}^{\EB}(\hat{\bm\eta}(\hat{\bm\theta}^{\ML})))]$ is denoted the excess variance caused by $\hat{\bm\eta}(\hat{\bm\theta}^{\ML})$. This is the additional cost due to the hyper-parameter estimation. It will vanish when we consider an EB estimator with a deterministic hyper-parameter, e.g., $\hat{\bm\theta}^{\TR}(\bm\eta)$ as discussed in \eqref{eq:2nd order limit of MSE difference}, or a generalized Bayes estimator $\hat{\bm\theta}^{\Bayes}$ in \eqref{eq:def of Bayes estimator}, which contains no hyper-parameter. One may believe that this component is nonnegative as it is associated with the variance incurred by estimating parameters, and this is true for the scaled EB hyper-parameter estimator, see Section \ref{subsec: XMSE of RLS using EB}. However, this is not necessarily the case for other hyper-parameter estimators, see Section \ref{subsec: XMSE of RLS using SURE and GCV}. 
\end{itemize}

The XMSE expression for a generalized Bayes estimator $\hat{\bm\theta}^{\Bayes}$ in \eqref{eq:def of Bayes estimator} is a direct consequence of Theorem \ref{thm:XMSE of empirical Bayes estimator}. Since this family of estimators is of independent interest, we state this result explicitly.


\begin{corollary}\label{corollary:XMSE of generalized Bayes estimator}
Under Assumptions \ref{asp:input and noise}-\ref{asp:limit of PPN}, $\XMSE(\hat{\bm\theta}^{\Bayes})$ is in the form of \eqref{eq:expression of XMSE} with $\Tr[\XVarHPE(\hat{\bm\theta}^{\Bayes})]=0$.
\end{corollary}

As for the regularized estimator $\hat{\bm\theta}^{\TR}(\hat{\bm\eta}(\hat{\bm\theta}^{\ML}))$ in \eqref{eq:RLS estimate}, we can also derive the expression of its XMSE, for which we need an additional assumption on the kernel matrix $\mathbf{P}(\bm\eta)$. 

\begin{assumption}\label{asp:kernel matrix}
The kernel matrix $\mathbf{P}(\bm\eta)$ is positive definite and twice continuously differentiable {with respect to} $\bm\eta$ in the interior of $\mathcal{D}_{\bm\eta}$.
\end{assumption}

\begin{corollary}\label{corollary:XSME of RLS estimator}
Let Assumptions \ref{asp:input and noise}-\ref{asp:kernel matrix} hold. The XMSE of $\hat{\bm\theta}^{\TR}(\hat{\bm\eta}(\hat{\bm\theta}^{\ML}))$ is in the form of \eqref{eq:expression of XMSE} with
\begin{align}
\label{eq:diff of b at theta0 for RLS}
&\bm{b}_{\bm\theta_{0},\star}(\bm\eta)=-\sigma^2\mathbf{\Sigma}^{-1}\mathbf{P}(\bm\eta)^{-1}\bm\theta_{0},\\
\label{eq:derivative of b wrt theta for RLS}
&\mathbf{b}_{\bm\theta_{0},\star}^{'}(\bm\eta)=-\sigma^2\mathbf{\Sigma}^{-1}\mathbf{P}(\bm\eta)^{-1},\\
\label{eq:derivative of b wrt eta for RLS}
&[\mathbf{b}_{\bm\eta,\star}^{'}(\bm\eta)]_{:,k}=-\sigma^2\mathbf{\Sigma}^{-1}\frac{\partial \mathbf{P}(\bm\eta)^{-1}}{\partial \eta_{k}}\bm\theta_{0}.
\end{align} 
\end{corollary}

By substituting \eqref{eq:derivative of b wrt theta for RLS} into \eqref{eq:def of XVar}, we {note} that for the regularized estimator $\hat{\bm\theta}^{\TR}(\hat{\bm\eta}(\hat{\bm\theta}^{\ML}))$, the excess variance caused by $\hat{\bm\theta}^{\TR}(\hat{\bm\eta}(\bm\theta_{0}))$ is negative, i.e.,
\begin{align}\label{eq:sign of XVar of RLS}
&\Tr[\XVar(\hat{\bm\theta}^{\TR}(\hat{\bm\eta}(\hat{\bm\theta}^{\ML})))]\nonumber\\
=&-2(\sigma^2)^2\Tr[\mathbf{\Sigma}^{-1}\mathbf{P}(\bm\eta_{\star}(\bm\theta_{0}))^{-1}\mathbf{\Sigma}^{-1}]{< 0}.
\end{align}
In fact, $\Tr[\XVar(\hat{\bm\theta}^{\TR}(\hat{\bm\eta}(\hat{\bm\theta}^{\ML})))]$ can also be expressed as
\begin{align*}
&\Tr[\XVar(\hat{\bm\theta}^{\TR}(\hat{\bm\eta}(\hat{\bm\theta}^{\ML})))]=\\
&\lim_{N\to\infty}N^2\{\Tr[\Var(\hat{\bm\theta}^{\TR}(\hat{\bm\eta}(\bm\theta_{0})))]-\Tr[\Var(\hat{\bm\theta}^{\ML})]\}{<0}.
\end{align*}
{This term is thus a manifestation of the rationale of regularization as mitigating the variance of the ML estimator.}

Moreover, if $\hat{\bm\eta}(\hat{\bm\theta}^{\ML})$ can be expressed as 
\begin{align}\label{eq:def of general hyper-parameter estimator}
\hat{\bm\eta}(\hat{\bm\theta}^{\ML})=\argmin_{\bm\eta\in\mathcal{D}_{\bm\eta}}\mathscr{F}_{N}(\hat{\bm\theta}^{\ML},\bm\eta),
\end{align}
and satisfies some additional assumptions, we can provide more detailed expressions of $\bm\eta_{\star}(\bm\theta_{0})$ and $\mathbf{D}_{\bm\eta,\star}^{'}(\bm\theta_{0})$.

\begin{assumption}\label{asp:compact set}
\begin{enumerate}
\item[]
\item The set $\mathcal{D}_{\bm\eta}$ is compact.
\item The hyper-parameter estimate $\hat{\bm\eta}(\hat{\bm\theta}^{\ML})$ satisfies the first-order optimality condition
\begin{align}\label{eq:first-order optimality condition}
\left.\frac{\partial \mathscr{F}_{N}(\hat{\bm\theta}^{\ML},\bm\eta)}{\partial \bm\eta}\right|_{\bm\eta=\hat{\bm\eta}(\hat{\bm\theta}^{\ML})}=\bm{0}.
\end{align}
\item There exists $\widetilde{\mathscr{F}}_{N}(\hat{\bm\theta}^{\ML},\bm\eta)=a^{\mathscr{F}}(\hat{\bm\theta}^{\ML})\mathscr{F}_{N}(\hat{\bm\theta}^{\ML},\bm\eta)+b^{\mathscr{F}}(\hat{\bm\theta}^{\ML})$
	with {$a^{\mathscr{F}}(\hat{\bm\theta}^{\ML})>0$} and $b^{\mathscr{F}}(\hat{\bm\theta}^{\ML})$ being independent of $\bm\eta$, such that $\widetilde{\mathscr{F}}_{N}(\bm\theta_{0},\bm\eta)$
 converges to a nonzero deterministic function $W(\bm\theta_{0},\bm\eta)$ uniformly for all $\bm\eta\in\mathcal{D}_{\bm\eta}$ as $N\to\infty$, i.e.,
	\begin{align}\label{eq:uniform as of F}
	\sup_{\bm\eta\in\mathcal{D}_{\bm\eta}}|\widetilde{\mathscr{F}}_{N}(\bm\theta_{0},\bm\eta)-W(\bm\theta_{0},\bm\eta)|\to 0,\ {\text{as}\ N\to\infty}.\quad 
	\end{align}
Moreover, when $\bm\eta=\argmin_{\bm\eta\in\mathcal{D}_{\bm\eta}}W(\bm\theta_{0},\bm\eta)$, we have ${\partial W(\bm\theta_{0},\bm\eta)}/{\partial \bm\eta}=\bm{0}$ and ${\partial^2 W(\bm\theta_{0},\bm\eta)}/{\partial \bm\eta \partial \bm\eta^{\top}}\succ 0$.
\end{enumerate}
\end{assumption}

\begin{assumption}\label{asp:convergence of cost function and derivatives}
{It holds that}
\begin{align}\label{eq:uniform convergence of A}
&\sup_{\bm\eta\in\mathcal{D}_{\bm\eta}}\left\|\frac{\partial^2 \widetilde{\mathscr{F}}_{N}(\bm\theta_{0},\bm\eta)}{\partial \bm\eta\partial \bm\eta^{\top}}-\mathbf{A}(\bm\theta_{0},\bm\eta) \right\|_{F}\to \mathbf{0},\\
\label{eq:uniform convergence of B}
&\sup_{\bm\eta\in\mathcal{D}_{\bm\eta}}\left\|\left.\frac{\partial^2 \widetilde{\mathscr{F}}_{N}(\hat{\bm\theta}^{\ML},\bm\eta)}{\partial \bm\eta\partial (\hat{\bm\theta}^{\ML})^{\top}}\right|_{\hat{\bm\theta}^{\ML}=\bm\theta_{0}}-\mathbf{B}(\bm\theta_{0},\bm\eta)\right\|_{F}\to 0,\quad
\end{align}
as $N\to\infty$, where $\mathbf{A}(\bm\theta_{0},\bm\eta)={\partial^2 W(\bm\theta_{0},\bm\eta)}/{\partial \bm\eta \partial \bm\eta^{\top}}$ and $\mathbf{B}(\bm\theta_{0},\bm\eta)={\partial^2 W(\bm\theta_{0},\bm\eta)}/{\partial \bm\eta \partial \bm\theta_{0}^{\top}}$.
\end{assumption}

\begin{remark}
Notice that the scaled EB, scaled SURE and scaled GCV hyper-parameter estimators introduced in Section \ref{subsec: FIR parameter estimators} all satisfy \eqref{eq:def of general hyper-parameter estimator} and \eqref{eq:uniform as of F}-\eqref{eq:uniform convergence of B}. For more details, we refer to Remark \ref{rmk:functions of ML} and Sections \ref{subsec:proof of expression of EB hyper-parameter estimator}-\ref{subsec:proof of XMSE of SY and GCV}.
\end{remark}

{We can now} provide more detailed expressions of $\bm\eta_{\star}(\bm\theta_{0})$ and $\mathbf{D}_{\bm\eta,\star}^{'}(\bm\theta_{0})$ defined in \eqref{eq:convergence of hyper-parameter estimator}-\eqref{eq:limit of derivative of hateta}.

\begin{corollary}\label{corollary:limit of derivative of eta}
Under Assumptions \ref{asp:compact set}-\ref{asp:convergence of cost function and derivatives}, we have
\begin{align}\label{eq:detailed expression of convergence of hyper-parameter estimator}
&\bm\eta_{\star}(\bm\theta_{0})=\argmin_{\bm\eta\in\mathcal{D}_{\bm\eta}}W(\bm\theta_{0},\bm\eta),\\
\label{eq:detailed expression of limit of derivative of hateta}
&\mathbf{D}_{\bm\eta,\star}^{'}(\bm\theta_{0})
=-\mathbf{A}(\bm\theta_{0},\bm\eta)^{-1}\mathbf{B}(\bm\theta_{0},\bm\eta)|_{\bm\eta=\bm\eta_{\star}(\bm\theta_{0})}.
\end{align}
\end{corollary}



{Next, based on Theorem \ref{thm:XMSE of empirical Bayes estimator} together with Corollaries \ref{corollary:XSME of RLS estimator}-\ref{corollary:limit of derivative of eta}, we will derive the explicit XMSE expressions for the regularized estimator with the hyper-parameter estimators detailed in Section \ref{subsec: FIR parameter estimators}. In fact, we only need to derive expressions of $\bm\eta_{\star}(\bm\theta_{0})$, $\mathbf{A}(\bm\theta_{0},\bm\eta)$ and $\mathbf{B}(\bm\theta_{0},\bm\eta)$ for the different estimators.}

\subsection{XMSE of the regularized estimator using the scaled EB hyper-parameter estimator}\label{subsec: XMSE of RLS using EB}

We start with the scaled EB hyper-parameter estimator.

\begin{corollary}\label{corollary:XMSE of EB}

Under Assumptions \ref{asp:input and noise}-\ref{asp:limit of PPN} and \ref{asp:kernel matrix}-\ref{asp:compact set}, for fixed $\alpha>0$, the XMSE of $\hat{\bm\theta}^{\TR}(\hat{\bm\eta}_{\EB,\alpha}))$ is in the form of \eqref{eq:expression of XMSE} with $\bm\eta_{\star}(\bm\theta_{0})=\bm\eta_{\tb\star,\alpha}$, $\mathbf{A}(\bm\theta_{0},\bm\eta)=\mathbf{A}_{\tb,\alpha}$ and $\mathbf{B}(\bm\theta_{0},\bm\eta)=\mathbf{B}_{\tb}$, where
\begin{subequations}\label{eq:expression for EB hyper-parameter estimator}
\begin{align}\label{eq:eta_b_star_alpha}
\bm\eta_{\tb\star,\alpha}=&\argmin_{\bm\eta\in\mathcal{D}_{\bm\eta}}W_{\tb,\alpha}(\bm\eta),\\
\label{eq:def of Wb_alpha}
W_{\tb,\alpha}(\bm\eta)=&\bm\theta_{0}^{\top}\mathbf{P}(\bm\eta)^{-1}\bm\theta_{0}+\alpha\log\det(\mathbf{P}(\bm\eta)),\\
\label{eq:def of Ab_alpha}
[\mathbf{A}_{\tb,\alpha}]_{k,l}=&(\bm\theta_{0})^{\top}\frac{\partial^2 \mathbf{P}(\eta)^{-1}}{\partial \eta_{k}\partial \eta_{l}}\bm\theta_{0}+\alpha\Tr\left[\mathbf{P}(\bm\eta)^{-1}\frac{\partial^2 \mathbf{P}(\bm\eta)}{\partial\eta_{k}\partial \eta_{l}} \right]\nonumber\\
&+\alpha\Tr\left[\frac{\partial \mathbf{P}(\bm\eta)^{-1}}{\partial \eta_{l}}\frac{\partial \mathbf{P}(\bm\eta)}{\partial \eta_{k}} \right],\\
\label{eq:def of Bb}
[\mathbf{B}_{\tb}]_{k,:}=&2(\bm\theta_{0})^{\top}\frac{\partial \mathbf{P}(\bm\eta)^{-1}}{\partial \eta_{k}}.
\end{align}
\end{subequations}
\end{corollary}

Notice that the excess variance caused by the scaled EB hyper-parameter estimator is nonnegative.


\begin{corollary}\label{corollary:sign of XVarHPE for eta_EB}
Under Assumptions \ref{asp:input and noise}-\ref{asp:limit of PPN} and \ref{asp:kernel matrix}-\ref{asp:compact set}, for $\hat{\bm\eta}_{\EB,\alpha}$, we have \eqref{eq:def of VarHPE part1} with $\mathbf{B}_{\tb}\mathbf{\Sigma}^{-1}=-2\mathbf{b}_{\bm\eta,\star}^{'}(\bm\eta)^{\top}$ and
\begin{align}\label{eq:XVarHPE for EB}
&\Tr[\XVarHPE(\hat{\bm\theta}^{\TR}(\hat{\bm\eta}_{\EB,\alpha}))]\\
=&4(\sigma^2)^2\Tr[\mathbf{\Sigma}^{-1}\mathbf{B}_{\tb}^{\top}\mathbf{A}_{\tb,\alpha}^{-1}\mathbf{B}_{\tb}\mathbf{\Sigma}^{-1}]|_{\bm\eta=\bm\eta_{\tb\star,\alpha}}
\geq 0,
\end{align}
{where the equality only holds if $\mathbf{B}_{\tb}|_{\bm\eta=\bm\eta_{\tb\star,\alpha}}=\mathbf{0}$.}
\end{corollary}

\subsection{XMSE of the regularized estimator with the scaled SURE and GCV hyper-parameter estimators}\label{subsec: XMSE of RLS using SURE and GCV}

In this subsection, we consider the scaled $\SURE_{y}$ and GCV hyper-parameter estimators. As {discussed} in Remark \ref{rmk:functions of ML}, the dominant terms of the cost functions of $\hat{\bm\eta}_{\Sy,\alpha}$ \eqref{eq:SUREy alpha hyper-parameter estimator} and $\hat{\bm\eta}_{\GCV,\alpha}$ \eqref{eq:GCV alpha hyper-parameter estimator} are the same. As a result, $\hat{\bm\eta}_{\Sy,\alpha}$ and $\hat{\bm\eta}_{\GCV,\alpha}$ share the same limit, which is similar to \cite[Theorem 1]{MCL18asy} and \cite[Theorem 1]{MCL2018gcv}. Consequently, $\hat{\bm\theta}^{\TR}(\hat{\bm\eta}_{\Sy,\alpha})$ and $\hat{\bm\theta}^{\TR}(\hat{\bm\eta}_{\GCV,\alpha})$ also share the same XMSE. The expressions of their corresponding $\bm\eta_{\star}(\bm\theta_{0})$, $\mathbf{A}(\bm\theta_{0},\bm\eta)$ and $\mathbf{B}(\bm\theta_{0},\bm\eta)$ are given in the following.

\begin{corollary}\label{corollary:XMSE of SUREy and GCV}
Under Assumptions \ref{asp:input and noise}-\ref{asp:limit of PPN} and \ref{asp:kernel matrix}-\ref{asp:compact set}, the XMSE of $\hat{\bm\theta}^{\TR}(\hat{\bm\eta}_{\Sy,\alpha})$ and $\hat{\bm\theta}^{\TR}(\hat{\bm\eta}_{\GCV,\alpha})$ is the same. It takes the form of \eqref{eq:expression of XMSE} with $\bm\eta_{\star}(\bm\theta_{0})=\bm\eta_{\ty\star,\alpha}$, $\mathbf{A}(\bm\theta_{0},\bm\eta)=\mathbf{A}_{\ty,\alpha}$ and $\mathbf{B}(\bm\theta_{0},\bm\eta)=\mathbf{B}_{\ty}$, where
\begin{subequations}
\begin{align}
\label{eq:opt hyperparameter of Wy_alpha}
\bm\eta_{\ty\star,\alpha}=&\argmin_{\bm\eta\in\mathcal{D}_{\bm\eta}}W_{\ty,\alpha}(\bm\eta),\\
\label{eq:def of Wy_alpha}
W_{\ty,\alpha}(\bm\eta)=&(\sigma^2)^2\bm\theta_{0}^{\top}\mathbf{P}(\bm\eta)^{-1}\mathbf{\Sigma}^{-1}\mathbf{P}(\bm\eta)^{-1}\bm{\theta}_{0}\nonumber\\
  &-2\alpha(\sigma^2)^2\Tr(\mathbf{\Sigma}^{-1}\mathbf{P}(\bm\eta)^{-1}),\\
\label{eq:def of Ay_alpha}
[\mathbf{A}_{\ty,\alpha}]_{k,l}
		=&2(\sigma^2)^2\left\{\bm\theta_{0}^{\top}\frac{\partial \mathbf{P}(\bm\eta)^{-1}}{\partial\eta_{l}}{\mathbf{\Sigma}}^{-1}\frac{\partial \mathbf{P}(\bm\eta)^{-1}}{\partial\eta_{k}}\bm\theta_{0}\right.\nonumber\\
		&+\bm\theta_{0}^{\top}\mathbf{P}(\bm\eta)^{-1}{\mathbf{\Sigma}}^{-1}\frac{\partial^2 \mathbf{P}(\bm\eta)^{-1}}{\partial\eta_{k}\partial\eta_{l}}\bm{\theta}_{0}\nonumber\\
		&-\alpha\left.\Tr\left[\mathbf{\Sigma}^{-1}\frac{\partial^2 \mathbf{P}(\bm\eta)^{-1}}{\partial\eta_{k}\partial\eta_{l}}\right]\right\},\\
        \label{eq:def of By}
  [\mathbf{B}_{\ty}]_{k,:}=&2({\sigma^2})^2\bm\theta_{0}^{\top}\mathbf{P}(\bm\eta)^{-1}\mathbf{\Sigma}^{-1}\frac{\partial \mathbf{P}(\bm\eta)^{-1}}{\partial\eta_{k}}\nonumber\\
&+2({\sigma^2})^2\bm\theta_{0}^{\top}\frac{\partial \mathbf{P}(\bm\eta)^{-1}}{\partial\eta_{k}}\mathbf{\Sigma}^{-1}\mathbf{P}(\bm\eta)^{-1}.
\end{align}
\end{subequations}
\end{corollary}

As it is generally acknowledged that there is a variance penalty associated with estimating model parameters \cite{Ljung1999}, one may believe that the excess variance caused by estimating $\bm\eta$ is positive. This actually holds under the conditions in the next corollary.



\begin{corollary}\label{corollary:sign of VarHPE of Sy GCV}
Under Assumptions \ref{asp:input and noise}-\ref{asp:limit of PPN} and \ref{asp:kernel matrix}-\ref{asp:compact set}, if one of the following conditions is satisfied,
\begin{enumerate}
\item $\mathbf{\Sigma}=\mathbf{I}_{n}$, $\mathbf{P}(\eta)=\eta \mathbf{K}$ with $\eta>0$ and fixed $\mathbf{K}\succ 0$,
\item $\mathbf{P}(\eta)=\eta \mathbf{I}_{n}$ with $\eta>0$,
\item $\mathbf{\Sigma}=\diag\{\bm{s}\}$ and $\mathbf{P}(\bm\eta)=\diag\{\bm\eta\}$ with $\bm{s},\bm{\eta}>0$,
\end{enumerate}
then $\Tr[\XVarHPE(\cdot)]> 0$ for $\hat{\bm\theta}^{\TR}(\hat{\bm\eta}_{\Sy,\alpha})$ or $\hat{\bm\theta}^{\TR}(\hat{\bm\eta}_{\GCV,\alpha})$. 
\end{corollary}

Therefore, under the conditions of Corollary \ref{corollary:sign of VarHPE of Sy GCV}, there is a penalty, as measured by the XMSE, for estimating hyper-parameters $\bm\eta$. However, surprisingly at least to the authors, there are settings where estimating $\bm\eta$ reduces the XMSE. We have the following example.


\begin{example}\label{example:negative VarHPE}
Consider the SS kernel \eqref{eq:SS kernel} with $\gamma=0.5$, $\mathbf{\Sigma}=\diag\{[10,500]\}$ and $\bm\theta_{0}=[2.53,1]^{\top}$. By utilizing \eqref{eq:XVarHPE of SY and GCV}, we have $\Tr[\XVarHPE(\cdot)]\approx -1.45\alpha(\sigma^2)^2<0$ for $\hat{\bm\theta}^{\TR}(\hat{\eta}_{\Sy,\alpha})$ and $\hat{\bm\theta}^{\TR}(\hat{\eta}_{\GCV,\alpha})$.
\end{example}

\section{A Finite Sample Improvement of the XMSE}\label{sec:accuracy of XMSE}

In Corollaries \ref{corollary:XMSE of EB} and \ref{corollary:XMSE of SUREy and GCV}, we have derived explicit XMSE expressions for $\hat{\bm\theta}^{\TR}(\hat{\bm\eta}_{\EB,\alpha})$, $\hat{\bm\theta}^{\TR}(\hat{\bm\eta}_{\Sy,\alpha})$ and $\hat{\bm\theta}^{\TR}(\hat{\bm\eta}_{\GCV,\alpha})$. Even if the XMSE is asymptotic in nature, we expect that $\XMSE(\hat{\bm\theta}^{\TR}(\hat{\bm\eta}(\hat{\bm\theta}^{\ML})))/N^2$ {provides a reliable approximation of} the sample difference between $\MSE(\hat{\bm\theta}^{\TR}(\hat{\bm\eta}(\hat{\bm\theta}^{\ML})))$ and $\MSE(\hat{\bm\theta}^{\ML})$ for finite, {but not too small,} sample sizes. 

To check the validity of using the XMSE for finite sample sizes, we first revisit the motivating example in Section \ref{subsec:motivating example}. We then calculate the difference: sample $\MSE(\cdot)$ $-$ sample $\MSE(\hat{\bm\theta}^{\ML})$, denoted by ``sample $\Delta\text{MSE}(\cdot)$'', and also calculate ${\XMSE(\cdot)}/{N^2}$ of the three regularized estimators in Table~\ref{table:XMSE and apx. XMSE in the motivating example}. We evaluate the accuracy of ${\XMSE(\cdot)}/{N^2}$ by
\begin{align*}
&{acc}(\XMSE(\cdot))=\\
&100\times \left(1-\frac{|{\XMSE(\cdot)}/{N^2}-\text{sample}\ \Delta\text{MSE}(\cdot)|}{|\text{sample}\ \Delta\text{MSE}(\cdot)|} \right).
\end{align*}
The larger $acc(\XMSE(\cdot))$, the more accurate $\XMSE(\cdot)$. When ${\XMSE(\cdot)}/{N^2}$ fits the sample MSE difference perfectly, we have $acc(\XMSE(\cdot))=100$. As displayed in Table~\ref{table:XMSE and apx. XMSE in the motivating example}, the accuracy of XMSE expressions in this example is far from satisfactory. Especially, when we consider $\hat{\bm\theta}^{\TR}(\hat{\eta}_{\Sy})$ or $\hat{\bm\theta}^{\TR}(\hat{\eta}_{\GCV})$, the corresponding XMSE is negative, indicating that $\hat{\bm\theta}^{\TR}(\hat{\eta}_{\Sy})$ and $\hat{\bm\theta}^{\TR}(\hat{\eta}_{\GCV})$ will perform better than $\hat{\bm\theta}^{\ML}$. It is contradictory to their positive sample MSE differences in Table~\ref{table:XMSE and apx. XMSE in the motivating example} and the result shown in Fig. \ref{fig:motivating example}. The inaccuracy of XMSE is due to the {rather small} sample size $N=50$. However, rather than increasing $N$, we will improve the accuracy {of the XMSE approximations} for smaller $N$.

\begin{table}[!htbp]
\centering
\caption{Sample $\Delta\text{MSE}(\cdot)$, $\frac{\XMSE(\cdot)}{N^2}$ and apx. $\frac{\XMSE(\cdot)}{N^2}$ for $\hat{\bm\theta}^{\TR}(\hat{\eta}_{\EB})$, $\hat{\bm\theta}^{\TR}(\hat{\eta}_{\Sy})$ and $\hat{\bm\theta}^{\TR}(\hat{\eta}_{\GCV})$ in the motivating example}
\label{table:XMSE and apx. XMSE in the motivating example}
\resizebox{\hsize}{!}{
\begin{tabular}{cccc}
\hline
& $\hat{\bm\theta}^{\TR}(\hat{\eta}_{\EB})$ & $\hat{\bm\theta}^{\TR}(\hat{\eta}_{\Sy})$ & $\hat{\bm\theta}^{\TR}(\hat{\eta}_{\GCV})$ \\
\hline
sample $\Delta\text{MSE}(\cdot)$ & $2.76\times 10^{-1}$ & $2.99\times 10^{-2}$ & $3.11\times 10^{-2}$ \\
\hline
${\XMSE(\cdot)}/{N^2}$ & $1.86\times 10^{-2}$ & \multicolumn{2}{c}{$-3.18\times 10^{-3}$}\\
$acc(\XMSE(\cdot))$ & $6.76$ & $-10.64$ & $-10.22$\\
\hline
apx. ${\XMSE(\cdot)}/{N^2}$ & $3.05\times 10^{-1}$ & \multicolumn{2}{c}{$4.51\times 10^{-2}$}\\
$acc(\text{apx.}\XMSE(\cdot))$ & $89.38$ & $49.26$ & $55.16$ \\
\hline
\end{tabular}}
\end{table}

Theorem \ref{thm:XMSE of empirical Bayes estimator} shows that the XMSE of $\hat{\bm\theta}^{\TR}(\hat{\bm\eta}(\hat{\bm\theta}^{\ML}))$ essentially depends on limits in \eqref{eq:limit of PPN}, \eqref{eq:uniform convergence of A}-\eqref{eq:uniform convergence of B}, \eqref{eq:convergence of hyper-parameter estimator}-\eqref{eq:limit of derivative of hateta} and
\eqref{eq:diff of b at theta0}-\eqref{eq:derivative of b wrt eta}. This suggests that the inaccuracy of XMSE in this example comes from that the limits are not sufficiently close to their corresponding finite sample expressions for $N=50$. To reduce this problem, we replace these limits in XMSE expressions by their corresponding finite sample expressions, e.g., $\mathbf{\Sigma}^{-1}$ in XMSE is replaced by $N(\mathbf{\Phi}^{\top}\mathbf{\Phi})^{-1}$. We call such approximations by ``apx. XMSE''. Returning to the motivating example, Table~\ref{table:XMSE and apx. XMSE in the motivating example} shows apx. ${\XMSE(\cdot)}/{N^2}$ and corresponding $acc(\cdot)$. We can observe that the values of apx. XMSE are much closer to the sample MSE differences and thus much more accurate than the exact XMSE. In particular, the sign of apx. ${\XMSE(\cdot)}/{N^2}$ of $\hat{\bm\theta}^{\TR}(\hat{\eta}_{\Sy})$ and $\hat{\bm\theta}^{\TR}(\hat{\eta}_{\GCV})$ is consistent with those of their sample MSE differences. 


\begin{remark}\label{rmk:convergence rates of EB and SUREy}
In Table~\ref{table:XMSE and apx. XMSE in the motivating example}, $acc(\text{apx}.\XMSE(\cdot))$ of $\hat{\bm\theta}^{\TR}(\hat{\eta}_{\Sy})$ is smaller than that of $\hat{\bm\theta}^{\TR}(\hat{\eta}_{\EB})$, which indicates that the apx. XMSE of $\hat{\bm\theta}^{\TR}(\hat{\eta}_{\Sy})$ is still less accurate than that of $\hat{\bm\theta}^{\TR}(\hat{\eta}_{\EB})$. This coincides with \cite[Theorem 2]{MCL18asy} {stating that} the convergence rate of $\hat{\eta}_{\Sy}$ is slower than that of $\hat{\eta}_{\EB}$. Accordingly, the accuracy of apx. XMSE of $\hat{\bm\theta}^{\TR}(\hat{\eta}_{\Sy})$ may require more samples than that of $\hat{\bm\theta}^{\TR}(\hat{\eta}_{\EB})$.
\end{remark}

To demonstrate the accuracy of XMSE and apx. XMSE, we {generate} more systems and display their corresponding apx. XMSE in the following example.

\begin{example}\label{example:accuracy of XMSE}

We randomly generate $20$ collections of $\bm\theta_{0}$ and $\{u(t),y(t)\}_{t=1}^{N}$ using the method in Section \ref{subsec:motivating example} with $n=20$ and $N=50,100,150,200$. Then we demonstrate the sample $\Delta\text{MSE}(\cdot)$ of $\hat{\bm\theta}^{\TR}(\hat{\eta}_{\EB})$ and $\hat{\bm\theta}^{\TR}(\hat{\eta}_{\Sy})$, and their corresponding ${\XMSE(\cdot)}/{N^2}$ and apx. ${\XMSE(\cdot)}/{N^2}$ in Fig. \ref{fig:XMSE with apx}. Note that for $\hat{\eta}_{\EB}$ and $\hat{\eta}_{\Sy}$, we still consider the SS kernel \eqref{eq:SS kernel} with $\gamma=0.95$. We can observe that the apx. ${\XMSE(\cdot)}/{N^2}$ is closer to the sample $\Delta\text{MSE}(\cdot)$, and thus more accurate than the exact ${\XMSE(\cdot)}/{N^2}$ for both $\hat{\bm\theta}^{\TR}(\hat{\eta}_{\EB})$ and $\hat{\bm\theta}^{\TR}(\hat{\eta}_{\Sy})$.


\begin{figure}[!htbp]
\begin{subfigure}[t]{0.24\textwidth}
\centering
\includegraphics[width=\textwidth]{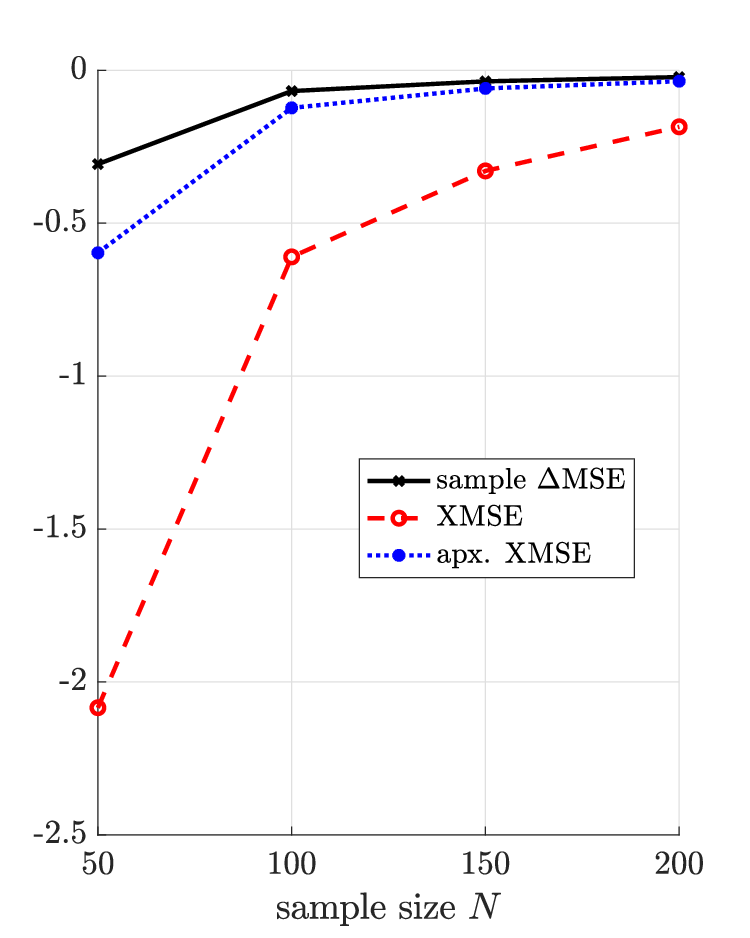}
\caption{}
\label{fig:apx XMSE EB}
\end{subfigure}%
\begin{subfigure}[t]{0.24\textwidth}
\centering
\includegraphics[width=\textwidth]{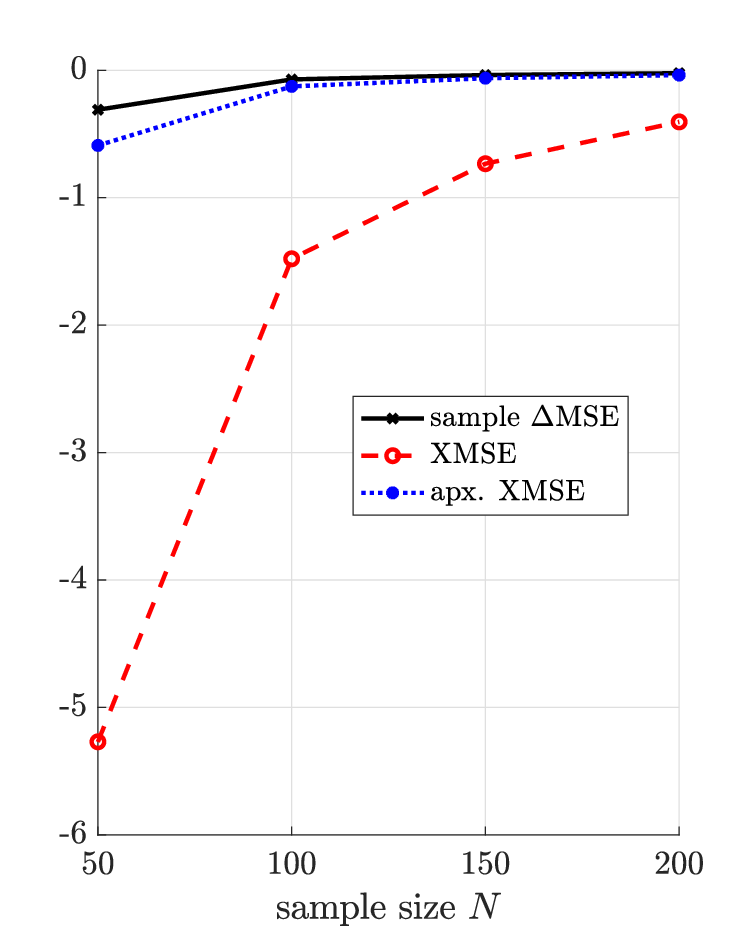}
\caption{}
\label{fig:apx XMSE SY}
\end{subfigure}
\caption{Sample $\Delta\text{MSE}(\cdot)$ (``sample $\Delta$MSE''),  ${\XMSE(\cdot)}/{N^2}$ (``XMSE''), and apx. ${\XMSE(\cdot)}/{N^2}$ (``apx. XMSE'') of $\hat{\bm\theta}^{\TR}(\hat{\eta}_{\EB})$ in Panel (a) and of $\hat{\bm\theta}^{\TR}(\hat{\eta}_{\Sy})$ in Panel (b).}
\label{fig:XMSE with apx}
\end{figure} 

\end{example}

\section{Numerical Simulation}\label{sec: numerical simulation}

In this section, we first generate systems with positive $\XMSE(\hat{\bm\theta}^{\TR}(\hat{\eta}_{\EB}))$ by using the XMSE expression given in Section \ref{sec:XMSE of empirical Bayes estimator}, and then study the XMSE components of the EB-based, SURE-based and GCV-based regularized estimators.

\subsection{Test systems and inputs}\label{subsec:test sys for n=20} 

We will generate $100$ collections of test systems and inputs with positive $\XMSE(\hat{\bm\theta}^{\TR}(\hat{\eta}_{\EB}))$. For each collection, we
\begin{enumerate}
\item generate a system with impulse response $\tilde{\bm\theta}_{0}\in\R^{n}$ and a collection of inputs $\{u(t)\}$ using the method in Section \ref{subsec:motivating example} with $n=20$ and $N=50$;
\item set $\sigma^2=1$, scale $\bm\theta_{0}=m\tilde{\bm\theta}_{0}$ such that the sample SNR is $10$ and select the SS kernel in \eqref{eq:SS kernel} with $\gamma=0.95$;
\item check whether or not $\XMSE(\hat{\bm\theta}^{\TR}(\hat{\eta}_{\EB}))>0.1$ is satisfied, if not then repeat the previous steps.
\end{enumerate}

\subsection{Simulation setup}\label{subsec:setup for positive XMSE}

For each collection of test systems and inputs, we conduct $100$ MC simulations to show the average performance of $\hat{\bm\theta}^{\ML}$, $\hat{\bm\theta}^{\TR}(\hat{\eta}_{\EB})$, $\hat{\bm\theta}^{\TR}(\hat{\eta}_{\Sy})$ and $\hat{\bm\theta}^{\TR}(\hat{\eta}_{\GCV})$. To evaluate the average performance of {an} estimator $\hat{\bm\theta}$, we consider two criteria: one is the sample $\MSE(\hat{\bm\theta})$ as defined after \eqref{eq:squared error of theta_hat}, the other one is the average $\FIT(\hat{\bm\theta})$ \cite{Ljung1995}, i.e., the sample mean of 
\begin{align}\label{eq: Fit}
\FIT(\hat{\bm\theta})=100\times \left(1-\frac{\|\hat{\bm\theta}-\bm\theta_{0}\|_{2}}{\|\bm\theta_{0}-\bar{\theta}_{0}\|_{2}} \right)
\end{align}
among $100$ MC simulations, where $\bar{\theta}_{0}=\sum_{k=1}^{n}[\bm\theta_{0}]_{k}/n$.
The better $\hat{\bm\theta}$ performs, the larger its average {FIT} and the smaller its sample MSE. The average {FIT} can be seen as a relative version of the sample MSE, but they may not have exact one-to-one correspondence. To better show the performance of estimators, we will provide both the average FIT and sample MSE. 

\subsection{Simulation results}\label{subsec:simulation results}

\begin{figure}[!htbp]
\centering
\includegraphics[width=0.9\linewidth]{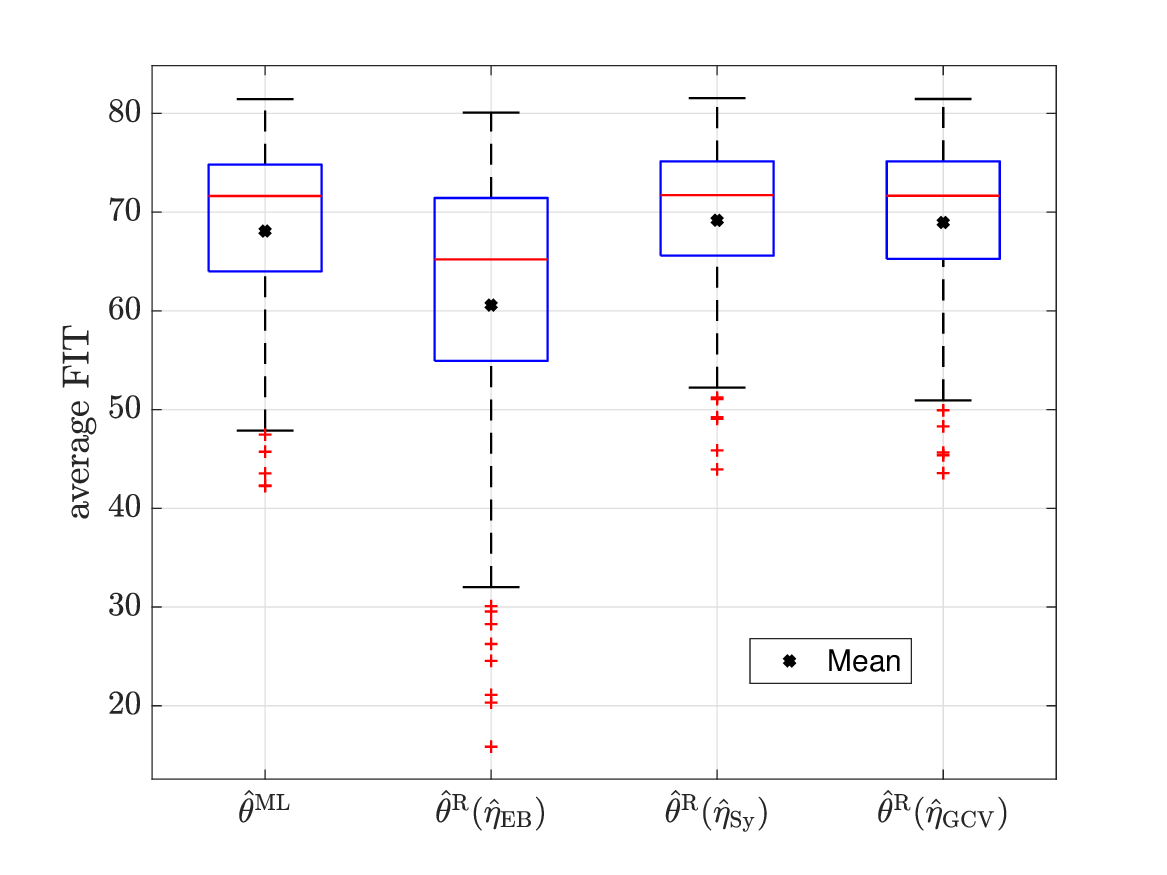}
\caption{Boxplots of average FIT of $\hat{\bm\theta}^{\ML}$, $\hat{\bm\theta}^{\TR}(\hat{\eta}_{\EB})$, $\hat{\bm\theta}^{\TR}(\hat{\eta}_{\Sy})$ and $\hat{\bm\theta}^{\TR}(\hat{\eta}_{\GCV})$ for systems with positive $\XMSE(\hat{\bm\theta}^{\TR}(\hat{\eta}_{\EB}))$.}
\label{fig:systems with positive XMSE KEEB}
\end{figure}

In Fig.~\ref{fig:systems with positive XMSE KEEB} and Table~\ref{table:average performances for systems with positive XMSE KEEB}, we demonstrate the average performance of $\hat{\bm\theta}^{\ML}$, $\hat{\bm\theta}^{\TR}(\hat{\eta}_{\EB})$, $\hat{\bm\theta}^{\TR}(\hat{\eta}_{\Sy})$ and $\hat{\bm\theta}^{\TR}(\hat{\eta}_{\GCV})$. We can observe that the average performance of $\hat{\bm\theta}^{\TR}(\hat{\eta}_{\EB})$ is worse than that of $\hat{\bm\theta}^{\ML}$, while the average performance of $\hat{\bm\theta}^{\TR}(\hat{\eta}_{\Sy})$ and $\hat{\bm\theta}^{\TR}(\hat{\eta}_{\GCV})$ is better. This is analogous to the motivating example in Section \ref{subsec:motivating example}. Moreover, among $100$ systems with positive $\XMSE(\hat{\bm\theta}^{\TR}(\hat{\eta}_{\EB}))$, there are $99$ systems where the sample $\MSE(\hat{\bm\theta}^{\TR}(\hat{\eta}_{\EB}))$ is larger than the sample $\MSE(\hat{\bm\theta}^{\ML})$ is $99$, which confirms the accuracy with the XMSE expression of $\hat{\bm\theta}^{\TR}(\hat{\eta}_{\EB})$. Also, there are $42$ and $45$ systems where the sample $\MSE(\hat{\bm\theta}^{\TR}(\hat{\eta}_{\Sy}))$ and $\MSE(\hat{\bm\theta}^{\TR}(\hat{\eta}_{\GCV}))$ are larger than the sample $\MSE(\hat{\bm\theta}^{\ML})$, respectively. This is consistent with the observation that the XMSE of $\hat{\bm\theta}^{\TR}(\hat{\eta}_{\Sy}))$ or $\hat{\bm\theta}^{\TR}(\hat{\eta}_{\GCV})$ for $46$ systems is positive.

\begin{table}[!htbp]
\centering
\caption{Means of average FIT and sample MSE of $\hat{\bm\theta}^{\ML}$, $\hat{\bm\theta}^{\TR}(\hat{\eta}_{\EB})$, $\hat{\bm\theta}^{\TR}(\hat{\eta}_{\Sy})$ and $\hat{\bm\theta}^{\TR}(\hat{\eta}_{\GCV})$ for systems with positive $\XMSE(\hat{\bm\theta}^{\TR}(\hat{\eta}_{\EB}))$.}
\label{table:average performances for systems with positive XMSE KEEB}
\resizebox{\hsize}{!}{
\begin{tabular}{ccccc}
\hline
& $\hat{\bm\theta}^{\ML}$ & $\hat{\bm\theta}^{\TR}(\hat{\eta}_{\EB})$ & $\hat{\bm\theta}^{\TR}(\hat{\eta}_{\Sy})$ & $\hat{\bm\theta}^{\TR}(\hat{\eta}_{\GCV})$ \\
\hline
average FIT & $68.23$ & $60.94$ & $69.36$ & $69.17$ \\
sample MSE & $9.01\times 10^{-1}$ & $1.37$ & $8.64\times 10^{-1}$ & $8.76\times 10^{-1}$ \\
\hline
\end{tabular}}
\end{table}

\begin{figure}[!htbp]
\centering
\includegraphics[width=\linewidth]{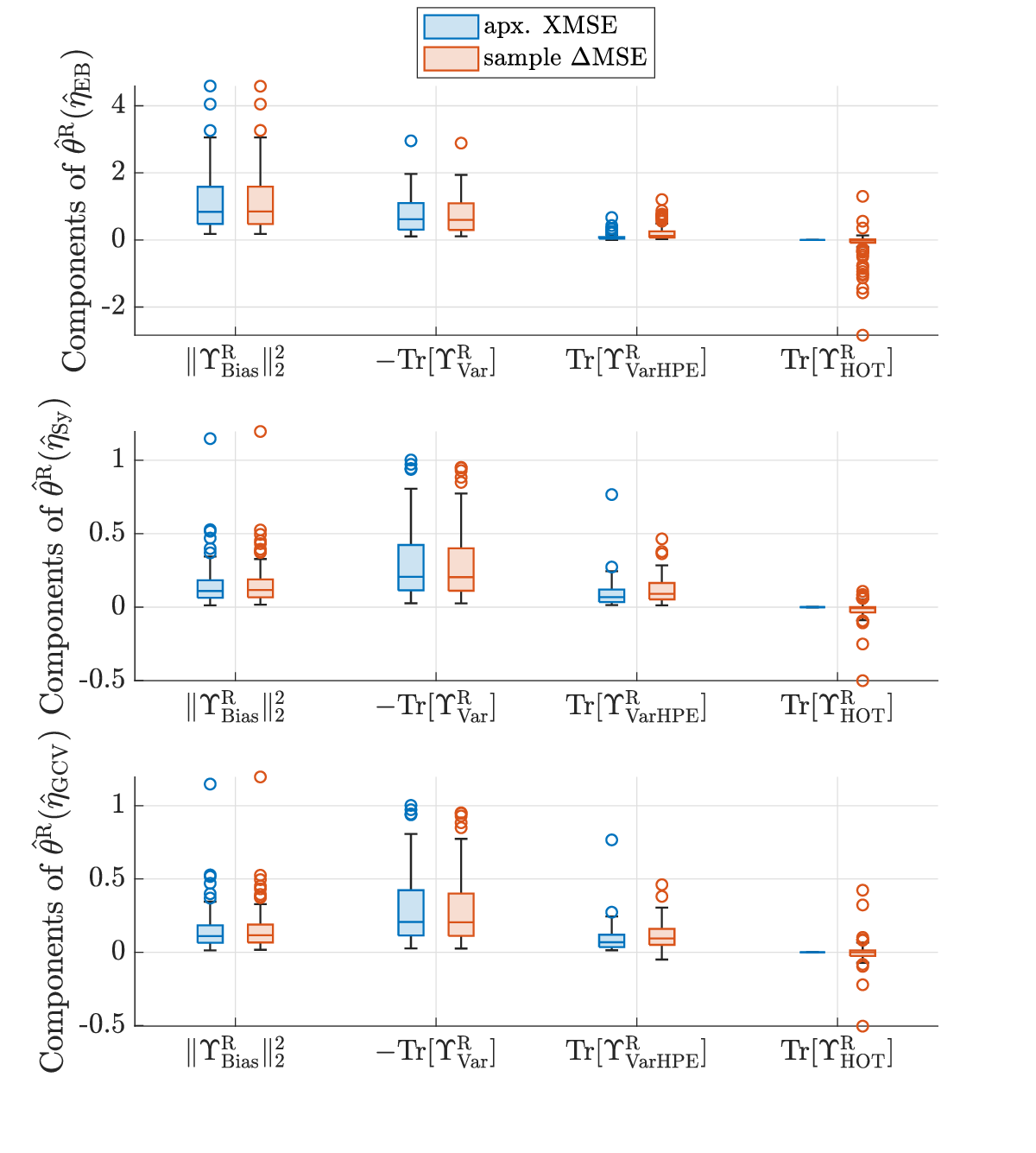}
\caption{Components of sample $\Delta\text{MSE}$ and apx. XMSE of $\hat{\bm\theta}^{\TR}(\hat{\eta}_{\EB})$, $\hat{\bm\theta}^{\TR}(\hat{\eta}_{\Sy})$ and $\hat{\bm\theta}^{\TR}(\hat{\eta}_{\GCV})$ for systems with positive $\XMSE(\hat{\bm\theta}^{\TR}(\hat{\eta}_{\EB}))$.}
\label{fig:XMSE components with positive XMSE KEEB}
\end{figure}

In Fig.~\ref{fig:XMSE components with positive XMSE KEEB}, we illustrate the components of sample $\Delta\text{MSE}(\cdot)$ and apx. $\XMSE(\cdot)/N^2$ of three regularized estimators. Their correspondence has been shown in Theorem \ref{thm:decomposition of MSE difference} and summarized in Table~\ref{table:components in apx. XMSE}. For regularized estimators, we use $\mathbf{\Upsilon}_{\Bias}^{\TR}$, $\mathbf{\Upsilon}_{\Var}^{\TR}$, $\mathbf{\Upsilon}_{\VarHPE}^{\TR}$ and $\mathbf{\Upsilon}_{\HOT}^{\TR}$ to denote components in Theorem \ref{thm:decomposition of MSE difference}, respectively. From Fig.~\ref{fig:XMSE components with positive XMSE KEEB}, we can notice that among three XMSE components, the apx. $\Tr[\XVarHPE(\cdot)]/{N^2}$ of all regularized estimators contributes the least, but the contributions of the other two components are different.
\begin{itemize}
\item[-] For $\hat{\bm\theta}^{\TR}(\hat{\eta}_{\EB})$, the apx. ${\|\XBias(\cdot)\|_{2}^2}/{N^2}>0$ contributes the most, leading to positive XMSE;
\item[-] for $\hat{\bm\theta}^{\TR}(\hat{\eta}_{\Sy})$ and $\hat{\bm\theta}^{\TR}(\hat{\eta}_{\GCV})$, the apx. $\Tr[\XVar(\cdot)]/N^2<0$ contributes the most, and thus, more than half of the systems still have negative XMSE. 
\end{itemize}
As shown in Table~\ref{table:XBias2, XVar and XVarHPE for motivating example}, the contributions of the XMSE components in the motivating example are similar. This implies that for either this experiment or the motivating example in Section \ref{subsec:motivating example}, the alignment between the true model parameter $\bm\theta_{0}$ and the kernel matrix $\mathbf{P}(\eta)=\eta\mathbf{K}$ with $[\mathbf{K}]_{k,l}=\kappa_{SS}(k,l;c=1,\gamma=0.95)$ defined in \eqref{eq:SS kernel} may not be that good. When this alignment is far from satisfactory, the XMSE of $\hat{\bm\theta}^{\TR}(\hat{\eta}_{\EB})$ can be positive and larger than that of $\hat{\bm\theta}^{\TR}(\hat{\eta}_{\Sy})$ or $\hat{\bm\theta}^{\TR}(\hat{\eta}_{\GCV})$.

\begin{table}[!htbp]
\caption{Components of sample $\Delta\text{MSE}(\hat{\bm\theta}^{\TR}(\hat{\bm\eta}(\hat{\bm\theta}^{\ML})))$ and apx. ${\XMSE(\hat{\bm\theta}^{\TR}(\hat{\bm\eta}(\hat{\bm\theta}^{\ML})))}/{N^2}$.}
\label{table:components in apx. XMSE}
\resizebox{0.9\hsize}{!}{
\begin{tabular}{l|l}
\hline
sample $\Delta$MSE & apx. XMSE\\
\hline
sample $\|\mathbf{\Upsilon}_{\Bias}^{\TR}\|_{2}^2$ & apx. $\|\XBias(\hat{\bm\theta}^{\TR}(\hat{\bm\eta}(\hat{\bm\theta}^{\ML})))\|_{2}^2/N^2$\\
sample $\Tr[\mathbf{\Upsilon}_{\Var}^{\TR}]$ & apx. $\Tr[\XVar(\hat{\bm\theta}^{\TR}(\hat{\bm\eta}(\hat{\bm\theta}^{\ML})))]/N^2$\\
sample $\Tr[\mathbf{\Upsilon}_{\VarHPE}^{\TR}]$ &
apx. $\Tr[\XVarHPE(\hat{\bm\theta}^{\TR}(\hat{\bm\eta}(\hat{\bm\theta}^{\ML})))]/N^2$\\
sample $\Tr[\mathbf{\Upsilon}_{\HOT}^{\TR}]$ & $0$\\
\hline
\end{tabular}}
\end{table}

\begin{table}[!htbp]
\caption{Components of sample $\Delta\text{MSE}$ and apx. $\XMSE$ of $\hat{\bm\theta}^{\TR}(\hat{\eta}_{\EB})$, $\hat{\bm\theta}^{\TR}(\hat{\eta}_{\Sy})$ and $\hat{\bm\theta}^{\TR}(\hat{\eta}_{\GCV})$ in the motivating example.}
\label{table:XBias2, XVar and XVarHPE for motivating example}
\centering
\resizebox{0.9\hsize}{!}{
\begin{tabular}{cccc}
\hline
& $\hat{\bm\theta}^{\TR}(\hat{\eta}_{\EB})$ & $\hat{\bm\theta}^{\TR}(\hat{\eta}_{\Sy})$ & $\hat{\bm\theta}^{\TR}(\hat{\eta}_{\GCV})$\\
\hline
 $\|\mathbf{\Upsilon}_{\Bias}^{\TR}\|_{2}^2$ & $5.05\times 10^{-1}$ &  \multicolumn{2}{c}{$5.86\times 10^{-2}$} \\
 $\|\XBias(\cdot)\|_{2}^2/N^2$ & $5.10\times 10^{-1}$ & \multicolumn{2}{c}{$5.93\times 10^{-2}$} \\
 \hline
 $\Tr[\mathbf{\Upsilon}_{\Var}^{\TR}]$ & $-2.37\times 10^{-1}$ & \multicolumn{2}{c}{$-7.30\times 10^{-2}$} \\
 $\Tr[\XVar(\cdot)]/N^2$ & $-2.46\times 10^{-1}$ & \multicolumn{2}{c}{$-7.59\times 10^{-2}$} \\
 \hline
 $\Tr[\mathbf{\Upsilon}_{\VarHPE}^{\TR}]$ & $5.92\times 10^{-2}$ & $4.97\times 10^{-2}$ & $5.00\times 10^{-2}$ \\
 $\Tr[\XVarHPE(\cdot)]/N^2$ & $4.11\times 10^{-2}$ & \multicolumn{2}{c}{$6.17\times 10^{-2}$} \\
 \hline
 $\Tr[\mathbf{\Upsilon}_{\HOT}^{\TR}]$ & $-4.63\times 10^{-2}$ & $-4.40\times 10^{-4}$ & $3.96\times 10^{-4}$ \\
\hline
\end{tabular}}
\end{table}

\begin{figure}[!htbp]
\centering
\includegraphics[width=\linewidth]{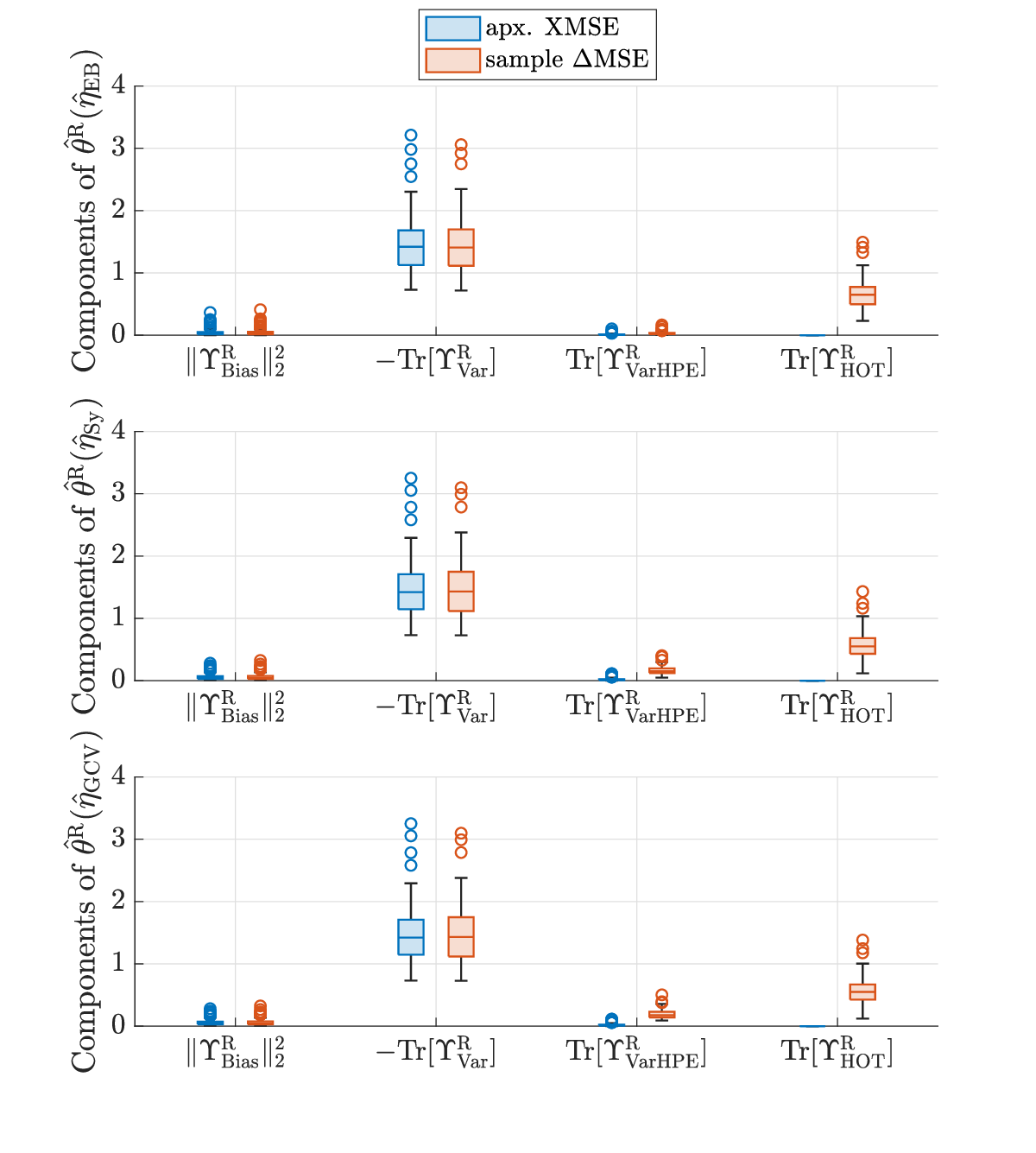}
\caption{Components of sample $\Delta\text{MSE}$ and apx. XMSE of $\hat{\bm\theta}^{\TR}(\hat{\eta}_{\EB})$, $\hat{\bm\theta}^{\TR}(\hat{\eta}_{\Sy})$ and $\hat{\bm\theta}^{\TR}(\hat{\eta}_{\GCV})$ for systems with proper alignment between $\bm\theta_{0}$ and $\mathbf{P}(\eta)=\eta\mathbf{K}$.}
\label{fig:XMSE components with proper alignment}
\end{figure}

To confirm our conjecture on the influence of alignment, we recall that the Gaussian weighting function is considered for the regularized estimator \eqref{eq:RLS estimate}. To achieve the proper alignment between $\bm\theta_{0}$ and $\mathbf{P}(\eta)=\eta\mathbf{K}$, we first generate $\{\tilde{\bm\theta}_{0,k}\}_{k=1}^{100}$ as $100$ independent realizations from $\mathcal{N}(\bm{0},\mathbf{K})$ with $[\mathbf{K}]_{k,l}=\kappa_{SS}(k,l;c=1,\gamma=0.95)$ and then follow the second step in Section \ref{subsec:test sys for n=20} to scale $\bm\theta_{0,k}=m_{k}\tilde{\bm\theta}_{0,k}$. From Fig.~\ref{fig:systems with positive XMSE KEEB} and Table~\ref{table:average performances for systems with proper alignment}, we can see that under the proper alignment, $\hat{\bm\theta}^{\TR}(\hat{\eta}_{\EB})$ performs best, and the performance of all three regularized estimators is noticeably better than $\hat{\bm\theta}^{\ML}$. We can also observe from Fig. \ref{fig:XMSE components with proper alignment} that, under the proper alignment, the apx. $\Tr[\XVar(\cdot)]/N^2<0$ contributes the most for $\hat{\bm\theta}^{\TR}(\hat{\eta}_{\EB})$, $\hat{\bm\theta}^{\TR}(\hat{\eta}_{\Sy})$ and $\hat{\bm\theta}^{\TR}(\hat{\eta}_{\GCV})$, leading to negative XMSE. 

\begin{remark}
It can be noticed from Tables~\ref{table:average performances for systems with positive XMSE KEEB} and \ref{table:average performances for systems with proper alignment} that the average performance of $\hat{\bm\theta}^{\TR}(\hat{\eta}_{\Sy})$ is slightly better than that of $\hat{\bm\theta}^{\TR}(\hat{\eta}_{\GCV})$. As numerically demonstrated in \cite{MCL2018gcv}, when the measurement noise variance $\sigma^2$ is unknown and needs to be estimated from the observed data, the GCV-based regularized estimator often slightly outperforms the SURE-based one. It is also mentioned in \cite{MCL2018gcv} that the possible reason is that the GCV-based regularized estimator can benefit from not requiring an estimate of $\sigma^2$. For the setting of this paper, we assume that $\sigma^2$ is known, and thus $\hat{\bm\theta}^{\TR}(\hat{\eta}_{\GCV})$ may diminish its advantage and perform slightly worse than $\hat{\bm\theta}^{\TR}(\hat{\eta}_{\Sy})$, as shown in Tables~\ref{table:average performances for systems with positive XMSE KEEB} and \ref{table:average performances for systems with proper alignment}.
\end{remark}

\begin{figure}[!htbp]
\centering
\includegraphics[width=0.9\linewidth]{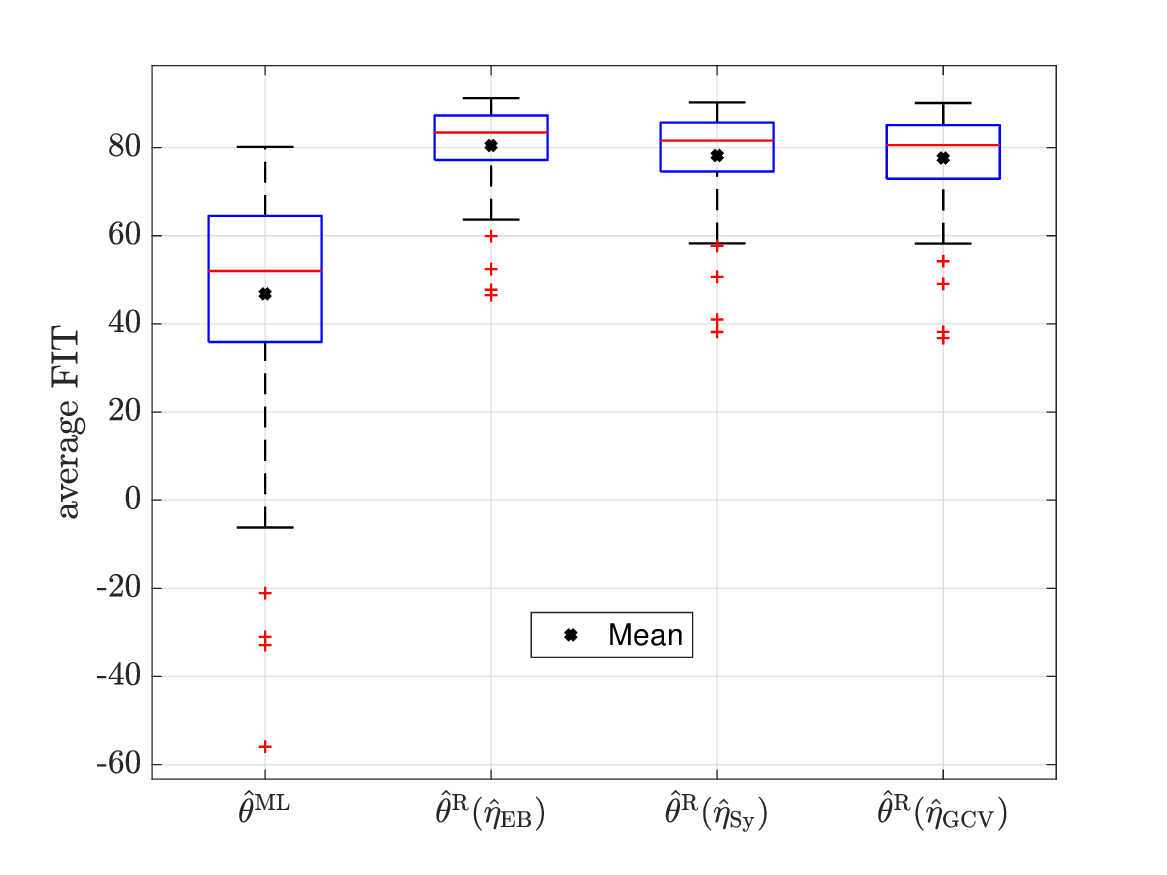}
\caption{Boxplots of average FIT of $\hat{\bm\theta}^{\ML}$, $\hat{\bm\theta}^{\TR}(\hat{\eta}_{\EB})$, $\hat{\bm\theta}^{\TR}(\hat{\eta}_{\Sy})$ and $\hat{\bm\theta}^{\TR}(\hat{\eta}_{\GCV})$ for systems with proper alignment between $\bm\theta_{0}$ and $\mathbf{P}(\eta)=\eta\mathbf{K}$.}
\label{fig:systems with proper alignment}
\end{figure}

\begin{remark}
Notice that in Theorem \ref{thm:XMSE of empirical Bayes estimator}, we have proved that the limit of $N^2\Tr[\mathbf{\Upsilon}_{\HOT}^{\TR}]$ is zero. However, $\Tr[\mathbf{\Upsilon}_{\HOT}^{\TR}]$ in Fig.~\ref{fig:XMSE components with positive XMSE KEEB}-\ref{fig:XMSE components with proper alignment} seems not always negligible compared with the other three components. We can approximate this term by deriving its higher order limit, i.e., $\lim_{N\to\infty}N^3\Tr[\mathbf{\Upsilon}_{\HOT}^{\TR}]$. Nevertheless, the corresponding expression will be more complicated and may make the analysis more challenging. Hence, we will not take this any further.
\end{remark}

As suggested by Fig.~\ref{fig:systems with positive XMSE KEEB}-\ref{fig:XMSE components with proper alignment} and Tables~\ref{table:average performances for systems with positive XMSE KEEB}-\ref{table:average performances for systems with proper alignment}, if we expect that the kernel-based regularized estimators can outperform the ML estimator, it is crucial to design a kernel in proper alignment with the true model parameters. 

\begin{table}[!htbp]
\centering
\caption{Means of average FIT and sample MSE of $\hat{\bm\theta}^{\ML}$, $\hat{\bm\theta}^{\TR}(\hat{\eta}_{\EB})$, $\hat{\bm\theta}^{\TR}(\hat{\eta}_{\Sy})$ and $\hat{\bm\theta}^{\TR}(\hat{\eta}_{\GCV})$ for systems with proper alignment between $\bm\theta_{0}$ and $\mathbf{P}(\eta)=\eta\mathbf{K}$.}
\label{table:average performances for systems with proper alignment}
\resizebox{\hsize}{!}{
\begin{tabular}{ccccc}
\hline
& $\hat{\bm\theta}^{\ML}$ & $\hat{\bm\theta}^{\TR}(\hat{\eta}_{\EB})$ & $\hat{\bm\theta}^{\TR}(\hat{\eta}_{\Sy})$ & $\hat{\bm\theta}^{\TR}(\hat{\eta}_{\GCV})$ \\
\hline
average FIT & $46.82$ & $80.48$ & $78.21$ & $77.64$ \\
sample MSE & $8.60\times 10^{-1}$ & $1.59\times 10^{-1}$ & $1.92\times 10^{-1}$ & $2.04\times 10^{-1}$ \\
\hline
\end{tabular}}
\end{table}

\section{Conclusion}\label{sec:conclusion}

This paper has addressed the problem of quantifying estimation accuracy for finite impulse response models using empirical Bayes (EB) estimators with measurements subject to white Gaussian noise. We have derived an explicit expression for the excess mean squared error (XMSE) of an EB estimator equipped with a general data-dependent hyper-parameter estimator. As specific examples, we have provided explicit XMSE expressions for kernel-based regularized estimators, which are EB estimators using a Gaussian weighting function with a zero mean and a covariance matrix parameterized by the scaled EB, Stein’s unbiased risk estimation in relation to output prediction ($\SURE_{y}$), and generalized cross-validation (GCV) hyper-parameter estimators. To improve the accuracy of these XMSE expressions for finite sample sizes, we have proposed a modification that has been validated in simulation studies for moderate sample sizes. We have applied these XMSE expressions to construct systems where the maximum likelihood (ML) estimator performs better than the regularized estimator. This sheds light on the influence of the alignment between model parameters and the chosen kernel on estimation accuracy. 

We believe the contributions in this paper can help us to enhance the understanding of kernel-based regularized estimators and design new estimators. 
\begin{itemize}
\item In \cite{JWH25cdc}, submitted to the $64$th IEEE Conference on Decision and Control (CDC), we focus on ridge regression, and provide a family of generalized Bayes estimators and a family of biased estimators in closed form. They have the same XMSE as the EB-based regularized estimator but are free of hyper-parameters. As shown numerically, these estimators exhibit performance comparable to the EB-based regularized estimator, while their implementation is computationally more efficient.
\item In a paper currently in preparation, we will develop more applications based on XMSE expressions in this paper.
\begin{itemize}
\item To deepen the understanding of kernel-based regularized estimators, we start by the observation that $\bm\theta_{0}^{\top}\mathbf{P}(\bm\eta)^{-1}\bm\theta_{0}$ in \eqref{eq:def of Wb_alpha}, $\bm\theta_{0}^{\top}\mathbf{P}(\bm\eta)^{-1}\mathbf{\Sigma}^{-1}\mathbf{P}(\bm\eta)^{-1}\bm{\theta}_{0}$ in \eqref{eq:def of Wy_alpha} and $\bm\theta_{0}^{\top}\mathbf{P}(\bm\eta)^{-1}\mathbf{\Sigma}^{-2}\mathbf{P}(\bm\eta)^{-1}\bm{\theta}_{0}$ in \eqref{eq:2nd order limit of MSE difference} will play important roles in XMSE expressions for regularized estimators, and could be used to quantify the alignment between the true model parameter $\bm\theta_{0}$ and the kernel matrix $\mathbf{P}(\bm\eta)$. Then, we look to further construct measures of such alignment and analyze their influences on the performance of regularized estimators. This will provide a more detailed explanation for the motivating example in Section \ref{subsec:motivating example} and support the alignment conjecture mentioned in Section \ref{subsec:simulation results}. Moreover, we analyze the influence of the number of hyper-parameters on the excess variances caused by different hyper-parameter estimators, i.e., on \eqref{eq:def of XVarHPE}.
\item To design new estimators with better performance, we utilize XMSE expressions presented in Theorem \ref{thm:XMSE of empirical Bayes estimator} and Corollary \ref{corollary:XMSE of generalized Bayes estimator} to construct a generalized Bayes estimator that outperforms the EB-based regularized estimator. Furthermore, we tune the parameter $\alpha$ in the scaled hyper-parameter estimators as introduced in Section \ref{subsec: FIR parameter estimators} based on XMSE expressions, thereby improving the performance of corresponding regularized estimators.
\end{itemize}
\end{itemize}




\def\thesectiondis{\thesection.}                   
\def\thesubsectiondis{\thesection.\arabic{subsection}.}          
\def\thesubsubsectiondis{\thesubsection.\arabic{subsubsection}.}

\setcounter{subsection}{0}

\renewcommand{\thesection}{A}
\setcounter{theorem}{0}
\renewcommand{\thelemma}{A.\arabic{lemma}}

\renewcommand{\theequation}{A.\arabic{equation}}
\setcounter{equation}{0}

\renewcommand{\thesubsection}{\thesection.\arabic{subsection}}

\section*{Appendix A}\label{sec:Appendix A}

Proofs of Theorems and Corollaries are included in Appendix A. Note that the proofs of Corollaries \ref{corollary:XMSE of generalized Bayes estimator} and \ref{corollary:sign of XVarHPE for eta_EB} are straightforward and thus omitted.

\subsection{Proof of Lemma \ref{lemma:limits of MSEs of ML and RLS estimators}}

Under Assumption \ref{asp:input and noise}, we have $\MSE(\hat{\bm\theta}^{\ML})=\sigma^2(\mathbf{\Phi}^{\top}\mathbf{\Phi})^{-1}$. Then we can utilize \eqref{eq:limit of PPN} to derive \eqref{eq:1st order limit of MSE of ML estimator}. For $\hat{\bm\theta}^{\TR}(\bm\eta)$, we recall \eqref{eq:reformulation of theta_RLS} to obtain $\hat{\bm\theta}^{\TR}(\bm\eta)-\hat{\bm\theta}^{\ML}=-\sigma^2(\mathbf{\Phi}^{\top}\mathbf{\Phi})^{-1}\mathbf{S}(\bm\eta)^{-1}\hat{\bm\theta}^{\ML}$.
Then we rewrite the MSE difference between $\hat{\bm\theta}^{\TR}(\bm\eta)$ and $\hat{\bm\theta}^{\ML}$,
\begin{align*}
&\MSE(\hat{\bm\theta}^{\TR}(\bm\eta))-\MSE(\hat{\bm\theta}^{\ML})\\
=&\E[\|\hat{\bm\theta}^{\TR}(\bm\eta)-\hat{\bm\theta}^{\ML}\|_{2}^2]
+2\E[(\hat{\bm\theta}^{\TR}(\bm\eta)-\hat{\bm\theta}^{\ML})(\hat{\bm\theta}^{\ML}-\bm\theta_{0})]\\
=&(\sigma^2)^2\Tr\{\mathbf{S}(\bm\eta)^{-1}(\mathbf{\Phi}^{\top}\mathbf{\Phi})^{-2}\mathbf{S}(\bm\eta)^{-1}[\bm\theta_{0}\bm\theta_{0}^{\top}+\sigma^2(\mathbf{\Phi}^{\top}\mathbf{\Phi})^{-1}]\}\\
&-2(\sigma^2)^2\Tr[(\mathbf{\Phi}^{\top}\mathbf{\Phi})^{-2}\mathbf{S}(\bm\eta)^{-1}].
\end{align*}
By utilizing \eqref{eq:limit of PPN}, we can obtain \eqref{eq:1st order limit of MSE of RLS estimator} and \eqref{eq:2nd order limit of MSE difference}.

\subsection{Proof of Remark \ref{rmk:functions of ML}}

To derive \eqref{eq:reformulation of theta_RLS}, we first rewrite  $(\mathbf{\Phi}^{\top}\mathbf{\Phi}+\sigma^2\mathbf{P}(\bm\eta)^{-1})^{-1}$ in \eqref{eq:RLS estimate} as $(\mathbf{\Phi}^{\top}\mathbf{\Phi})^{-1}-\sigma^2(\mathbf{\Phi}^{\top}\mathbf{\Phi})^{-1}\mathbf{S}(\bm\eta)^{-1}(\mathbf{\Phi}^{\top}\mathbf{\Phi})^{-1}$ using the Woodbury matrix identity, and then substitute \eqref{eq:ML estimate}. The reformulation of $\mathscr{F}_{\EB,\alpha}(\bm\eta)$ is based on \cite[(A.19),(B.11)]{MCL18asy}; the reformulation of $\mathscr{F}_{\Sy,\alpha}(\bm\eta)$ is based on \cite[(A.15)]{MCL18asy} with $[\mathbf{\Phi}^{\top}\mathbf{\Phi}+\sigma^2\mathbf{P}(\bm\eta)^{-1}]^{-1}\mathbf{P}(\bm\eta)^{-1}=(\mathbf{\Phi}^{\top}\mathbf{\Phi})^{-1}\mathbf{S}(\bm\eta)^{-1}$; the reformulation of $\mathscr{F}_{\GCV,\alpha}(\bm\eta)$ is based on the proof of \cite[Theorem 1]{MCL2018gcv}.

\subsection{Proof of Theorem \ref{thm:decomposition of MSE difference}}

We first rewrite the MSE difference between $\hat{\bm\theta}^{\EB}(\hat{\bm\eta}(\hat{\bm\theta}^{\ML}))$ and $\hat{\bm\theta}^{\ML}$ as follows,
\begin{align*}
&\MSE(\hat{\bm\theta}^{\EB}(\hat{\bm\eta}(\hat{\bm\theta}^{\ML})))-\MSE(\hat{\bm\theta}^{\ML})\\
=&\|\E(\hat{\bm\theta}^{\EB}(\hat{\bm\eta}(\hat{\bm\theta}^{\ML})))-\bm\theta_{0}\|_{2}^2\\
&+\Tr[\Var(\hat{\bm\theta}^{\EB}(\hat{\bm\eta}(\hat{\bm\theta}^{\ML})))]
-\Tr[\Var(\hat{\bm\theta}^{\ML})],
\end{align*}
where we utilize $\E(\hat{\bm\theta}^{\ML})=\bm\theta_{0}$. Then we plug 
\begin{align*}
\hat{\bm\theta}^{\EB}(\hat{\bm\eta}(\hat{\bm\theta}^{\ML}))=&[\hat{\bm\theta}^{\EB}(\hat{\bm\eta}(\hat{\bm\theta}^{\ML}))-\hat{\bm\theta}^{\EB}(\hat{\bm\eta}(\bm\theta_{0}))]\\
&+[\hat{\bm\theta}^{\EB}(\hat{\bm\eta}(\bm\theta_{0}))-\hat{\bm\theta}^{\ML}]+\hat{\bm\theta}^{\ML}
\end{align*}
into $\Var(\hat{\bm\theta}^{\EB}(\hat{\bm\eta}(\hat{\bm\theta}^{\ML})))$ to obtain
\begin{align*}
&\Var(\hat{\bm\theta}^{\EB}(\hat{\bm\eta}(\hat{\bm\theta}^{\ML})))-\Var(\hat{\bm\theta}^{\ML})\\
=&\mathbf{\Upsilon}_{\VarHPE}+\mathbf{\Upsilon}_{\Var}+\Var(\hat{\bm\theta}^{\EB}(\hat{\bm\eta}(\hat{\bm\theta}^{\ML}))-\hat{\bm\theta}^{\ML}).
\end{align*}
We now complete the proof of Theorem \ref{thm:decomposition of MSE difference}.

\subsection{Proof of Theorem \ref{thm:XMSE of empirical Bayes estimator}}\label{subsec: proof of XMSE of empirical Bayes estimator}

We will derive limits of $N\mathbf{\Upsilon}_{\Bias}$, $N^2\mathbf{\Upsilon}_{\Var}$, $N^2\mathbf{\Upsilon}_{\VarHPE}$ and $N^2\mathbf{\Upsilon}_{\HOT}$ as defined in Theorem \ref{thm:decomposition of MSE difference}, respectively. 

To derive the limits of $N\mathbf{\Upsilon}_{\Bias}$ and $N^2\mathbf{\Upsilon}_{\Var}$, we apply the second-order Taylor expansion\footnote{Note that for random sequences $\xi_{N}\in\C$ and $a_{N}\in\C$, if $\xi_{N}/a_{N}$ converges in probability to $0$ as $N\to\infty$, we denote it as $\xi_{N}=o_{p}(a_{N})$.} of $\hat{\bm\theta}^{\EB}(\hat{\bm\eta}(\bm\theta_{0}))-\hat{\bm\theta}^{\ML}$ at $\hat{\bm\theta}^{\ML}=\bm\theta_{0}$ as follows,
\begin{align}\label{eq:diff thetaR_theta0 and theta_ML}
&\hat{\bm\theta}^{\EB}(\hat{\bm\eta}(\bm\theta_{0}))-\hat{\bm\theta}^{\ML}
=\hat{\bm\theta}^{\EB}(\hat{\bm\eta}(\bm\theta_{0}))|_{\hat{\bm\theta}^{\ML}=\bm\theta_{0}}-\bm\theta_{0}\nonumber\\
&\quad+\left(\left.\frac{\partial \hat{\bm\theta}^{\EB}(\hat{\bm\eta}(\bm\theta_{0}))}{\partial \hat{\bm\theta}^{\ML}}\right|_{\hat{\bm\theta}^{\ML}=\bm\theta_{0}}-\mathbf{I}_{n}\right)(\hat{\bm\theta}^{\ML}-\bm\theta_{0})\nonumber\\
&\quad+o_{p}(\|\hat{\bm\theta}^{\ML}-\bm\theta_{0}\|_{2}^2).
\end{align}
In the following, we will substitute $\hat{\bm\theta}^{\EB}(\bm\eta)$ for $\hat{\bm\theta}^{\EB}(\hat{\bm\eta}(\bm\theta_{0}))$ for simplicity. 

Then our next step is to derive the limits of $\hat{\bm\theta}^{\EB}(\bm\eta)-\bm\theta_{0}$ and ${\partial \hat{\bm\theta}^{\EB}(\bm\eta)}/{\partial \hat{\bm\theta}^{\ML}}-\mathbf{I}_{n}$ given ${\hat{\bm\theta}^{\ML}=\bm\theta_{0}}$. Recall the definition of $\hat{\bm\theta}^{\EB}(\bm\eta)$ in \eqref{eq: empirical Bayes estimator}. By applying the second-order Taylor expansion to $\log p(\bm{Y}|\bm{\theta})$ at $\bm\theta=\hat{\bm\theta}^{\ML}$, we have\footnote{Note that for sequences $\xi_{N}\in\C$ and $a_{N}\in\C$, if $\lim_{N\to\infty}\xi_{N}/a_{N}=0$, we denote it as $\xi_{N}=o(a_{N})$.} 
\begin{align*}
\log p(\bm{Y}|\bm{\theta})=&\log p(\bm{Y}|\hat{\bm{\theta}}^{\ML})\\
&+(1/\sigma^2)(\bm{Y}-\mathbf{\Phi}\hat{\bm\theta}^{\ML})^{\top}\mathbf{\Phi}(\bm\theta-\hat{\bm\theta}^{\ML})\\
&-[1/(2\sigma^2)](\bm\theta-\hat{\bm\theta}^{\ML})^{\top}(\mathbf{\Phi}^{\top}\mathbf{\Phi})(\bm\theta-\hat{\bm\theta}^{\ML})\\
&+o(\|(\mathbf{\Phi}^{\top}\mathbf{\Phi})^{-1}\|_{F}).
\end{align*}
Notice that $(\bm{Y}-\mathbf{\Phi}\hat{\bm\theta}^{\ML})^{\top}\mathbf{\Phi}(\bm\theta-\hat{\bm\theta}^{\ML})=0$, which can be derived from $\mathbf{\Phi}^{\top}\mathbf{\Phi}\hat{\bm\theta}^{\ML}=\mathbf{\Phi}^{\top}\bm{Y}$. It follows that $p(\bm{Y}|\bm\theta)=p(\bm{Y}|\hat{\bm{\theta}}^{\ML})\exp[-(1/2\sigma^2)(\bm\theta-\hat{\bm\theta}^{\ML})^{\top}(\mathbf{\Phi}^{\top}\mathbf{\Phi})(\bm\theta-\hat{\bm\theta}^{\ML})]+o(\|(\mathbf{\Phi}^{\top}\mathbf{\Phi})^{-1}\|_{F})$.
We then have
\begin{align}\noeqref{eq: Taylor expansion of empirical Bayes estimator}
\hat{\bm\theta}^{\EB}({\bm\eta})=&\frac{\int \bm\theta\pi(\bm\theta|{\bm\eta})p(\bm\theta|\hat{\bm\theta}^{\ML})d\bm\theta}{\int \pi(\bm\theta|{\bm\eta})p(\bm\theta|\hat{\bm\theta}^{\ML})d\bm\theta}+o(\|(\mathbf{\Phi}^{\top}\mathbf{\Phi})^{-1}\|_{F}),\nonumber\\
\label{eq: Taylor expansion of empirical Bayes estimator}
=&\frac{\E_{\bm\theta|\hat{\bm\theta}^{\ML}}[\bm\theta \pi(\bm\theta|\bm\eta)]}{\E_{\bm\theta|\hat{\bm\theta}^{\ML}}[\pi(\bm\theta|{\bm\eta})]}
+o(\|(\mathbf{\Phi}^{\top}\mathbf{\Phi})^{-1}\|_{F}),
\end{align}
where the expectation $\E_{\bm\theta|\hat{\bm\theta}^{\ML}}$ is with respect to $p(\bm\theta|\hat{\bm\theta}^{\ML})=\mathcal{N}(\hat{\bm\theta}^{\ML},\sigma^2(\mathbf{\Phi}^{\top}\mathbf{\Phi})^{-1})$. Next, we apply the second-order Taylor expansions to $\pi(\bm\theta|\bm\eta)$ and $\bm\theta\pi(\bm\theta|\bm\eta)$ at $\bm\theta=\hat{\bm\theta}^{\ML}$, and calculate their expectations with respect to $p(\bm\theta|\hat{\bm\theta}^{\ML})$ to obtain
\begin{align}\label{eq:Taylor expansion of expectation of pi}
\E_{\bm\theta|\hat{\bm\theta}^{\ML}}[\pi(\bm\theta|\bm\eta)]=&\pi(\hat{\bm\theta}^{\ML}|\bm\eta)\nonumber\\
&+\frac{\sigma^2}{2}\Tr\left[\left.\frac{\partial^2 \pi(\bm\theta|\bm\eta)}{\partial \bm\theta\partial\bm\theta^{\top}}\right|_{\bm\theta=\hat{\bm\theta}^{\ML}}(\mathbf{\Phi}^{\top}\mathbf{\Phi})^{-1}\right]\nonumber\\
&+o(\|(\mathbf{\Phi}^{\top}\mathbf{\Phi})^{-1}\|_{F}),\\
\label{eq:Taylor expansion of expectation of theta_pi}
\E_{\bm\theta|\hat{\bm\theta}^{\ML}}[\bm\theta\pi(\bm\theta|\bm\eta)]=&\hat{\bm\theta}^{\ML}\pi(\hat{\bm\theta}^{\ML}|\bm\eta)\nonumber\\
+&\sigma^2(\mathbf{\Phi}^{\top}\mathbf{\Phi})^{-1}\left.\frac{\partial \pi(\bm\theta|\bm\eta)}{\partial \bm\theta}\right|_{\bm\theta=\hat{\bm\theta}^{\ML}}\nonumber\\
+&\frac{\sigma^2}{2}\hat{\bm\theta}^{\ML}\Tr\left[\left.\frac{\partial^2 \pi(\bm\theta|\bm\eta)}{\partial \bm\theta\partial\bm\theta^{\top}}\right|_{\bm\theta=\hat{\bm\theta}^{\ML}}(\mathbf{\Phi}^{\top}\mathbf{\Phi})^{-1}\right]\nonumber\\
+&o(\|(\mathbf{\Phi}^{\top}\mathbf{\Phi})^{-1}\|_{F}).
\end{align}
It follows that
\begin{align*}
\hat{\bm\theta}^{\EB}(\bm\eta)
=&\frac{\E_{\bm\theta|\hat{\bm\theta}^{\ML}}[\bm\theta\pi(\bm\theta|\bm\eta)]/\pi(\hat{\bm\theta}^{\ML}|\bm\eta)}{1+[\E_{\bm\theta|\hat{\bm\theta}^{\ML}}[\pi(\bm\theta|\bm\eta)]/\pi(\hat{\bm\theta}^{\ML}|\bm\eta)-1]}\\
&+o(\|(\mathbf{\Phi}^{\top}\mathbf{\Phi})^{-1}\|_{F})\\
=&\frac{\E_{\bm\theta|\hat{\bm\theta}^{\ML}}[\bm\theta\pi(\bm\theta|\bm\eta)]}{\pi(\hat{\bm\theta}^{\ML}|\bm\eta)}
\left\{1-\left[\frac{\E_{\bm\theta|\hat{\bm\theta}^{\ML}}[\pi(\bm\theta|\bm\eta)]}{\pi(\hat{\bm\theta}^{\ML}|\bm\eta)}-1\right]\right\}\\
&+o(\|(\mathbf{\Phi}^{\top}\mathbf{\Phi})^{-1}\|_{F})\\
=&\hat{\bm\theta}^{\ML}+\sigma^2(\mathbf{\Phi}^{\top}\mathbf{\Phi})^{-1}\left.\frac{\partial \log(\pi(\bm\theta|\bm\eta))}{\partial \bm\theta}\right|_{\bm\theta=\hat{\bm\theta}^{\ML}}\\
&+o(\|(\mathbf{\Phi}^{\top}\mathbf{\Phi})^{-1}\|_{F}),
\end{align*}
where for the second step, we apply the first-order Taylor expansion of $1/\{1+[\E_{\bm\theta|\hat{\bm\theta}^{\ML}}[\pi(\bm\theta|\bm\eta)]/\pi(\hat{\bm\theta}^{\ML}|\bm\eta)-1]\}$ at $\E_{\bm\theta|\hat{\bm\theta}^{\ML}}[\pi(\bm\theta|\bm\eta)]/\pi(\hat{\bm\theta}^{\ML}|\bm\eta)-1=0$; for the last step, we utilize \eqref{eq:Taylor expansion of expectation of pi}-\eqref{eq:Taylor expansion of expectation of theta_pi} and 
\begin{align*}
\frac{\partial \log(\pi(\bm\theta|\bm\eta))}{\partial \bm\theta}=&\frac{1}{\pi(\bm\theta|\bm\eta)}\frac{\partial \pi(\bm\theta|\bm\eta)}{\partial \bm\theta}.
\end{align*}
Then we have
\begin{align}\label{eq: Taylor expansion of EB estimator}
\hat{\bm\theta}^{\EB}(\bm\eta)|_{\hat{\bm\theta}^{\ML}=\bm\theta_{0}}-\bm\theta_{0}=&\sigma^2(\mathbf{\Phi}^{\top}\mathbf{\Phi})^{-1}\left.\frac{\partial \log(\pi(\bm\theta|\bm\eta))}{\partial \bm\theta}\right|_{\bm\theta=\bm\theta_{0}}\nonumber\\
&+o(\|(\mathbf{\Phi}^{\top}\mathbf{\Phi})^{-1}\|_{F}),\\
\label{eq: Taylor expansion of EB estimator derivative}
\left.\frac{\partial \hat{\bm\theta}^{\EB}(\bm\eta)}{\partial \hat{\bm\theta}^{\ML}}\right|_{\hat{\bm\theta}^{\ML}=\bm\theta_{0}}-\mathbf{I}_{n}=&\sigma^2(\mathbf{\Phi}^{\top}\mathbf{\Phi})^{-1}\left.\frac{\partial^2 \log(\pi(\bm\theta|\bm\eta))}{\partial \bm\theta \partial \bm\theta^{\top}}\right|_{\bm\theta=\bm\theta_{0}}\nonumber\\
&+o(\|(\mathbf{\Phi}^{\top}\mathbf{\Phi})^{-1}\|_{F}).
\end{align}
By substituting \eqref{eq: Taylor expansion of EB estimator}-\eqref{eq: Taylor expansion of EB estimator derivative} into \eqref{eq:diff thetaR_theta0 and theta_ML}, and applying \eqref{eq:limit of PPN} and \eqref{eq:convergence of hyper-parameter estimator}, we obtain \eqref{eq:def of XBias}-\eqref{eq:def of XVar part1} and \eqref{eq:diff of b at theta0}-\eqref{eq:derivative of b wrt theta}.

To derive the limit of $N^2\mathbf{\Upsilon}_{\VarHPE}$, we consider the first-order Taylor expansion of $\hat{\bm\theta}^{\EB}(\hat{\bm\eta}(\hat{\bm\theta}^{\ML}))-\hat{\bm\theta}^{\EB}(\hat{\bm\eta}(\bm\theta_{0}))$ at $\hat{\bm\theta}^{\ML}=\bm\theta_{0}$ as follows,
\begin{align}
&\hat{\bm\theta}^{\EB}(\hat{\bm\eta}(\hat{\bm\theta}^{\ML}))-\hat{\bm\theta}^{\EB}(\hat{\bm\eta}(\bm\theta_{0}))\\
=&\left.\frac{\partial \hat{\bm\theta}^{\EB}(\hat{\bm\eta}(\hat{\bm\theta}^{\ML}))}{\partial \hat{\bm\eta}(\hat{\bm\theta}^{\ML})}\frac{\partial \hat{\bm\eta}(\hat{\bm\theta}^{\ML})}{\partial \hat{\bm\theta}^{\ML}}\right|_{\hat{\bm\theta}^{\ML}=\bm\theta_{0}}(\hat{\bm\theta}^{\ML}-\bm\theta_{0})\\
\label{eq:Taylor expansion of theta_etaML}
&+o_{p}(\|\hat{\bm\theta}^{\ML}-\bm\theta_{0}\|_{2}),
\end{align}
where
\begin{align}
&\left.\frac{\partial \hat{\bm\theta}^{\EB}(\hat{\bm\eta}(\hat{\bm\theta}^{\ML}))}{\partial \hat{\bm\theta}^{\ML}}\right|_{\hat{\bm\theta}^{\ML}=\bm\theta_{0}}
=\nonumber\\
&\left.\left[\frac{\partial \hat{\bm\theta}^{\EB}(\bm\eta)}{\partial \hat{\bm\theta}^{\ML}}
+\frac{\partial \hat{\bm\theta}^{\EB}(\bm\eta)}{\partial \bm\eta}\frac{\partial \hat{\bm\eta}(\hat{\bm\theta}^{\ML})}{\partial \hat{\bm\theta}^{\ML}}\right]\right|_{\bm\eta=\hat{\bm\eta}(\bm\theta_{0}),\hat{\bm\theta}^{\ML}=\bm\theta_{0}},\nonumber\\
\label{eq:derivative of theta_R wrt eta}
&\left.\frac{\partial \hat{\bm\theta}^{\EB}(\bm\eta)}{\partial \bm\eta}\right|_{\hat{\bm\theta}^{\ML}=\bm\theta_{0}}=\sigma^2(\mathbf{\Phi}^{\top}\mathbf{\Phi})^{-1}\left.\frac{\partial^2 \log(\pi(\bm\theta|\bm\eta))}{\partial \bm\theta \partial \bm\eta^{\top}}\right|_{\bm\theta=\bm\theta_{0}}\nonumber\\
&\qquad+o(\|(\mathbf{\Phi}^{\top}\mathbf{\Phi})^{-1}\|_{F}).
\end{align}
Note that \eqref{eq:derivative of theta_R wrt eta} is derived from \eqref{eq: Taylor expansion of EB estimator}. By using \eqref{eq:limit of PPN} and \eqref{eq:convergence of hyper-parameter estimator}-\eqref{eq:limit of derivative of hateta}, we obtain \eqref{eq:def of XVarHPE}-\eqref{eq:def of VarHPE part1} and \eqref{eq:derivative of b wrt eta}.

To derive the limit of $N^2\mathbf{\Upsilon}_{\HOT}$, we first rewrite $\mathbf{\Upsilon}_{\HOT}$ as $\mathbf{\Upsilon}_{\HOT}
=\E(\mathcal{D}_{\bm\eta}\mathcal{D}_{\bm\eta}^{\top})+\E(\mathcal{D}_{\bm\eta}\mathcal{D}_{\bm\theta}^{\top})+\E(\mathcal{D}_{\bm\theta}\mathcal{D}_{\bm\eta}^{\top})+\Var(\mathcal{D}_{\bm\theta})$ with $\mathcal{D}_{\bm\eta}=\hat{\bm\theta}^{\EB}(\hat{\bm\eta}(\hat{\bm\theta}^{\ML}))-\hat{\bm\theta}^{\EB}(\hat{\bm\eta}(\bm\theta_{0}))$ and $\mathcal{D}_{\bm\theta}=\hat{\bm\theta}^{\EB}(\hat{\bm\eta}(\bm\theta_{0}))-\hat{\bm\theta}^{\ML}$, where we use $\E(\hat{\bm\theta}^{\ML})=\bm\theta_{0}$ and $\hat{\bm\theta}^{\EB}(\hat{\bm\eta}(\hat{\bm\theta}^{\ML}))-\hat{\bm\theta}^{\ML}=\mathcal{D}_{\bm\eta}+\mathcal{D}_{\bm\theta}$.
Therefore, we apply \eqref{eq:diff thetaR_theta0 and theta_ML} and \eqref{eq:Taylor expansion of theta_etaML} together with \eqref{eq:limit of PPN} and \eqref{eq:limit of derivative of hateta} to obtain $\lim_{N\to\infty}N^2\mathbf{\Upsilon}_{\HOT}=0$.

\subsection{Proof of Corollary \ref{corollary:XSME of RLS estimator}}

To derive \eqref{eq:diff of b at theta0 for RLS}-\eqref{eq:derivative of b wrt eta for RLS}, we utilize the equivalent form of $\hat{\bm\theta}^{\TR}(\bm\eta)$ in \eqref{eq:reformulation of theta_RLS} and calculate the limits of $N[\hat{\bm\theta}^{\TR}(\bm\eta)|_{\hat{\bm\theta}^{\ML}=\bm\theta_{0}}-\bm\theta_{0}]$,
\begin{align*}
N\left[\left.\frac{\partial \hat{\bm\theta}^{\TR}(\bm\eta)}{\partial \hat{\bm\theta}^{\ML}}\right|_{\hat{\bm\theta}^{\ML}=\bm\theta_{0}}-\mathbf{I}_{n}\right]\ \text{and}\ 
N\left.\frac{\partial \hat{\bm\theta}^{\TR}(\bm\eta)}{\partial \bm\eta}\right|_{\hat{\bm\theta}^{\ML}=\bm\theta_{0}}.
\end{align*}

\subsection{Proof of Corollary \ref{corollary:limit of derivative of eta}}

As mentioned in Assumption \ref{asp:compact set}, both $a^{\mathscr{F}}(\hat{\bm\theta}^{\ML})$ and $b^{\mathscr{F}}(\hat{\bm\theta}^{\ML})$ are {independent} of $\bm\eta$. It leads to $\hat{\bm\eta}(\hat{\bm\theta}^{\ML})=\argmin_{\bm\eta\in\mathcal{D}_{\bm\eta}}{\mathscr{F}_{N}}(\hat{\bm\theta}^{\ML},\bm\eta)=\argmin_{\bm\eta\in\mathcal{D}_{\bm\eta}}\widetilde{\mathscr{F}}_{N}(\hat{\bm\theta}^{\ML},\bm\eta)$.
{We can then} apply \cite[Theorem 8.2]{Ljung1999} together with $\hat{\bm\eta}(\bm\theta_{0})=\argmin_{\bm\eta\in\mathcal{D}_{\bm\eta}}\widetilde{\mathscr{F}}_{N}(\bm\theta_{0},\bm\eta)$ and \eqref{eq:uniform as of F} to derive \eqref{eq:detailed expression of convergence of hyper-parameter estimator}.
Note that $\hat{\bm\eta}(\hat{\bm\theta}^{\ML})$ satisfies the first-order optimality condition \eqref{eq:first-order optimality condition}, leading to ${\partial \widetilde{\mathscr{F}}_{N}(\hat{\bm\theta}^{\ML},\bm\eta)}/{\partial \bm\eta}|_{\bm\eta=\hat{\bm\eta}(\hat{\bm\theta}^{\ML})}=\bm{0}$.
Take its full derivative with respect to $\hat{\bm\theta}^{\ML}$:
\begin{align*}
&\left.\frac{\partial^2 \widetilde{\mathscr{F}}_{N}(\hat{\bm\theta}^{\ML},\bm\eta)}{\partial \bm\eta\partial (\hat{\bm\theta}^{\ML})^{\top}}\right|_{\bm\eta=\hat{\bm\eta}(\hat{\bm\theta}^{\ML})}\nonumber\\
&+\left.\frac{\partial^2 \widetilde{\mathscr{F}}_{N}(\hat{\bm\theta}^{\ML},\bm\eta)}{\partial \bm\eta\partial \bm\eta^{\top}}\right|_{\bm\eta=\hat{\bm\eta}(\hat{\bm\theta}^{\ML})}\frac{\partial \hat{\bm\eta}(\hat{\bm\theta}^{\ML})}{\partial \hat{\bm\theta}^{\ML}}=\mathbf{0}.
\end{align*}
Thus, ${\partial \hat{\bm\eta}(\hat{\bm\theta}^{\ML})}/{\partial \hat{\bm\theta}^{\ML}}|_{\hat{\bm\theta}^{\ML}=\bm\theta_{0}}$ can be rewritten as
\begin{align*}
-\left(\frac{\partial^2 \widetilde{\mathscr{F}}_{N}(\bm\theta_{0},\bm\eta)}{\partial \bm\eta\partial \bm\eta^{\top}}\right)^{-1}
\left.\frac{\partial^2 \widetilde{\mathscr{F}}_{N}(\hat{\bm\theta}^{\ML},\bm\eta)}{\partial \bm\eta\partial (\hat{\bm\theta}^{\ML})^{\top}}\right|_{\bm\eta=\hat{\bm\eta}(\bm\theta_{0}),\hat{\bm\theta}^{\ML}=\bm\theta_{0}}.
\end{align*}
By using \cite[Theorem 8.2]{Ljung1999} and \eqref{eq:uniform convergence of A}-\eqref{eq:uniform convergence of B}, we obtain \eqref{eq:detailed expression of limit of derivative of hateta}.

\subsection{Proof of Corollary \ref{corollary:XMSE of EB}}\label{subsec:proof of expression of EB hyper-parameter estimator}

To derive the limit of $\hat{\bm\eta}_{\EB,\alpha}|_{\hat{\bm\theta}^{\ML}=\bm\theta_{0}}$, we first recall the reformulation of $\mathscr{F}_{\EB,\alpha}(\bm\eta)$ in Remark \ref{rmk:functions of ML} to obtain $\hat{\bm\eta}_{\EB,\alpha}=\argmin_{\bm\eta\in\mathcal{D}_{\bm\eta}}\mathscr{F}_{\tb,\alpha}(\hat{\bm\theta}^{\ML},\bm\eta)$. Then we prove the uniform convergence of $\mathscr{F}_{\tb,\alpha}(\bm\theta_{0},\bm\eta)$ to $W_{\tb,\alpha}(\bm\eta)$ for $\bm\eta\in\mathcal{D}_{\bm\eta}$, where $W_{\tb,\alpha}(\bm\eta)$ is defined in \eqref{eq:def of Wb_alpha}. For the difference between $\mathscr{F}_{\tb,\alpha}(\bm\theta_{0},\bm\eta)$ and $W_{\tb,\alpha}(\bm\eta)$, we utilize the norm inequalities in \cite[Chapter 10.4]{PP2012} to derive
\begin{align*}
&\sup_{\bm\eta\in\mathcal{D}_{\bm\eta}}|\mathscr{F}_{\tb,\alpha}(\bm\theta_{0},\bm\eta)-W_{\tb,\alpha}(\bm\eta)|\\
\leq&\|\bm\theta_{0}\|_{2}^2\sup_{\bm\eta\in\mathcal{D}_{\bm\eta}}\|\mathbf{S}(\bm\eta)^{-1}-\mathbf{P}(\bm\eta)^{-1}\|_{F}\\
&+\alpha\sup_{\bm\eta\in\mathcal{D}_{\bm\eta}}|\log\det(\mathbf{S}(\bm\eta))-\log\det(\mathbf{P}(\bm\eta))|.
\end{align*}
We only need to derive the uniform convergence of $\mathbf{S}(\bm\eta)^{-1}$ to $\mathbf{P}(\bm\eta)^{-1}$, and $\log\det(\mathbf{S}(\bm\eta))$ to $\log\det(\mathbf{P}(\bm\eta))$, respectively. Due to the compactness of $\mathcal{D}_{\bm\eta}$, we know
\begin{align}
\label{eq:boundedness of S and P norm}
&\sup_{\bm\eta\in\mathcal{D}_{\bm\eta}}\|\mathbf{S}(\bm\eta)^{-1}\|_{F}\leq\sup_{\bm\eta\in\mathcal{D}_{\bm\eta}}\|\mathbf{P}(\bm\eta)^{-1}\|_{F}\ \text{is}\ \text{bounded}.\quad\ 
\end{align}
Note $\mathbf{S}(\bm\eta)^{-1}-\mathbf{P}(\bm\eta)^{-1}=-\sigma^2\mathbf{S}(\bm\eta)^{-1}(\mathbf{\Phi}^{\top}\mathbf{\Phi})^{-1}\mathbf{P}(\bm\eta)^{-1}$. Then we have 
\begin{align}\label{eq:uniform convergence of S_inv}
&\sup_{\bm\eta\in\mathcal{D}_{\bm\eta}}\|\mathbf{S}(\bm\eta)^{-1}-\mathbf{P}(\bm\eta)^{-1}\|_{F}\leq\nonumber\\
& \sigma^2\|(\mathbf{\Phi}^{\top}\mathbf{\Phi})^{-1}\|_{F}\sup_{\bm\eta\in\mathcal{D}_{\bm\eta}}\|\mathbf{P}(\bm\eta)^{-1}\|_{F}^2
\to 0,
\end{align}
as $N\to\infty$. For the uniform convergence of $\log\det(\mathbf{S}(\bm\eta))$, we rewrite the difference between $\log\det(\mathbf{S}(\bm\eta))$ and $\log\det(\mathbf{P}(\bm\eta))$ as $\log\det(\mathbf{P}(\bm\eta)^{-1/2}\mathbf{S}(\bm\eta)\mathbf{P}(\bm\eta)^{-1/2})$.
Then we utilize \eqref{eq:limit of PPN}, the norm inequalities in \cite[Chapter 10.4]{PP2012}, \eqref{eq:boundedness of S and P norm} and \cite[Lemma B.16]{JCL2021} to obtain $\sup_{\bm\eta\in\mathcal{D}_{\bm\eta}}|\log\det(\mathbf{S}(\bm\eta))-\log\det(\mathbf{P}(\bm\eta))|\to 0$ as $N\to\infty$. Therefore, we can prove $\sup_{\bm\eta\in\mathcal{D}_{\bm\eta}}|\mathscr{F}_{\tb,\alpha}(\bm\theta_{0},\bm\eta)-W_{\tb,\alpha}(\bm\eta)|\to 0$ as $N\to\infty$. As a result of \cite[Theorem 8.2]{Ljung1999}, we have $\lim_{N\to\infty}\hat{\bm\eta}_{\EB,\alpha}|_{\hat{\bm\theta}^{\ML}=\bm\theta_{0}}
=\bm\eta_{\tb\star,\alpha}$.

The expression of $\mathbf{A}_{\tb,\alpha}$ and the uniform convergence in \eqref{eq:uniform convergence of A} can be similarly proved like \cite[Proof of Theorem 1]{JCML21eb}.
%
%
To derive the expression of $\mathbf{B}_{\tb}$, we first calculate $\mathbf{B}_{\tb}={\partial^2 W_{\tb,\alpha}(\bm\eta)}/{\partial \bm\eta\partial \bm\theta_{0}^{\top}}$ in \eqref{eq:def of Bb}
and then prove the uniform convergence of ${\partial^2 \mathscr{F}_{\tb,\alpha}(\bm\eta)}/{\partial \bm\eta\partial (\hat{\bm\theta}^{\ML})^{\top}}|_{\hat{\bm\theta}^{\ML}=\bm\theta_{0}}$. We have
\begin{align*}
&\sup_{\bm\eta\in\mathcal{D}_{\bm\eta}}\|{\partial^2 \mathscr{F}_{\tb,\alpha}(\bm\eta)}/{\partial \bm\eta\partial (\hat{\bm\theta}^{\ML})^{\top}}|_{\hat{\bm\theta}^{\ML}=\bm\theta_{0}}-\mathbf{B}_{\tb}\|_{F}\\
\leq&\sum_{k}\sup_{\bm\eta\in\mathcal{D}_{\bm\eta}}\|{\partial^2 \mathscr{F}_{\tb,\alpha}(\bm\eta)}/{\partial \eta_{k}\partial (\hat{\bm\theta}^{\ML})^{\top}}|_{\hat{\bm\theta}^{\ML}=\bm\theta_{0}}-[\mathbf{B}_{\tb}]_{k,:}\|_{2}\\
\leq& \sum_{k}2\|\hat{\bm\theta}^{\ML}-\bm\theta_{0}\|_{2}\sup_{\bm\eta\in\mathcal{D}_{\bm\eta}}\|{\partial \mathbf{S}(\bm\eta)^{-1}}/{\partial \eta_{k}}\|_{F}\\
&+2\|\bm\theta_{0}\|_{2} \sup_{\bm\eta\in\mathcal{D}_{\bm\eta}}\|{\partial \mathbf{S}(\bm\eta)^{-1}}/{\partial \eta_{k}}-{\partial \mathbf{P}(\bm\eta)^{-1}}/{\partial \eta_{k}}\|_{F}\to 0,
\end{align*} 
where ${\partial^2 \mathscr{F}_{\tb,\alpha}(\bm\eta)}/{\partial \eta_{k}\partial (\hat{\bm\theta}^{\ML})^{\top}}=2(\hat{\bm\theta}^{\ML})^{\top}[{\partial \mathbf{S}(\bm\eta)^{-1}}/{\partial \eta_{k}}]$. For the last step, we use \eqref{eq:boundedness of S and P norm}, \eqref{eq:limit of PPN}, \cite[(59), Chapter 10.4]{PP2012},
\begin{align}\label{eq:uniform convergence of derivative of S}
&\sup_{\bm\eta\in\mathcal{D}_{\bm\eta}}\|{\partial \mathbf{S}(\bm\eta)^{-1}}/{\partial \eta_{k}}-{\partial \mathbf{P}(\bm\eta)^{-1}}/{\partial \eta_{k}}\|_{F}\nonumber\\
\leq & 2\sup_{\bm\eta\in\mathcal{D}_{\bm\eta}}\|\mathbf{S}(\bm\eta)^{-1}-\mathbf{P}(\bm\eta)^{-1}\|_{F}\sup_{\bm\eta\in\mathcal{D}_{\bm\eta}}\|{\partial \mathbf{P}(\bm\eta)}/{\partial \eta_{k}}\|_{F}\nonumber\\
&\times \sup_{\bm\eta\in\mathcal{D}_{\bm\eta}}\|\mathbf{P}(\bm\eta)^{-1}\|_{F}\to 0,\ \text{as}\ N\to\infty,
\end{align}
and the boundedness of $\sup_{\bm\eta\in\mathcal{D}_{\bm\eta}}\left\|{\partial \mathbf{P}(\bm\eta)}/{\partial \eta_{k}}\right\|_{F}$, which is due to the compactness of $\mathcal{D}_{\bm\eta}$.

\subsection{Proof of Corollary \ref{corollary:XMSE of SUREy and GCV}}\label{subsec:proof of XMSE of SY and GCV}

To derive the limit of $\hat{\bm\eta}_{\Sy,\alpha}|_{\hat{\bm\theta}^{\ML}=\bm\theta_{0}}$, we recall Remark \ref{rmk:functions of ML} and know that $\hat{\bm\eta}_{\Sy,\alpha}=\argmin_{\bm\eta\in\mathcal{D}_{\bm\eta}}N\mathscr{F}_{\ty,\alpha}(\hat{\bm\theta}^{\ML},\bm\eta)$. Then we prove the uniform convergence of $N\mathscr{F}_{\ty,\alpha}(\bm\theta_{0},\bm\eta)$ to $W_{\ty,\alpha}(\bm\eta)$ defined in \eqref{eq:def of Wy_alpha}. Analogously to the proof in Section \ref{subsec:proof of expression of EB hyper-parameter estimator}, we can prove the uniform convergence of $N\mathscr{F}_{\ty,\alpha}(\bm\theta_{0},\bm\eta)$ to $W_{\ty,\alpha}(\bm\eta)$ based on the convergence of $\sup_{\bm\eta\in\mathcal{D}_{\bm\eta}}\|\mathbf{S}(\bm\eta)^{-1}-\mathbf{P}(\bm\eta)^{-1}\|_{F}$ and $\|N(\mathbf{\Phi}^{\top}\mathbf{\Phi})^{-1}-\mathbf{\Sigma}^{-1}\|_{F}$ to zero in \eqref{eq:uniform convergence of S_inv} and \eqref{eq:limit of PPN}. Using \cite[Theorem 8.2]{Ljung1999}, we have $\lim_{N\to\infty}\hat{\bm\eta}_{\Sy,\alpha}|_{\hat{\bm\theta}^{\ML}=\bm\theta_{0}}
=\bm\eta_{\ty\star,\alpha}$.

To derive the expression of $\mathbf{A}_{\ty,\alpha}$, we first calculate $\mathbf{A}_{\ty,\alpha}={\partial^2 W_{\ty,\alpha}(\bm\eta)}/{\partial \bm\eta\partial\bm\eta^{\top}}$ to obtain \eqref{eq:def of Ay_alpha}. Then we prove the uniform convergence of ${\partial^2 N\mathscr{F}_{\ty,\alpha}(\bm\theta_{0},\bm\eta)}/{\partial \bm\eta\partial\bm\eta^{\top}}$ to $\mathbf{A}_{\ty,\alpha}$, which is essentially the uniform convergence of their entries:
\begin{align*}
&\frac{\partial^2 N{\mathscr{F}}_{\ty,\alpha}(\bm\theta_{0},\bm\eta)}{\partial \eta_{k}\partial \eta_{l}}
=2({\sigma^2})^2\bm\theta_{0}^{\top}\mathbf{S}(\bm\eta)^{-1}N(\mathbf{\Phi}^{\top}\mathbf{\Phi})^{-1}\frac{\partial^2 \mathbf{S}(\bm\eta)^{-1}}{\partial \eta_{k}\partial\eta_{l}}\bm\theta_{0}\\
&\qquad\qquad+2({\sigma^2})^2\bm\theta_{0}^{\top}\frac{\mathbf{S}(\bm\eta)^{-1}}{\partial \eta_{l}}N(\mathbf{\Phi}^{\top}\mathbf{\Phi})^{-1}\frac{\partial \mathbf{S}(\bm\eta)^{-1}}{\partial \eta_{k}}\bm\theta_{0}\\
&\qquad\qquad-2\alpha({\sigma^2})^2\Tr\left[N(\mathbf{\Phi}^{\top}\mathbf{\Phi})^{-1}\frac{\partial^2 \mathbf{S}(\bm\eta)^{-1}}{\partial\eta_{k}\partial\eta_{l}} \right]
\end{align*}
and \eqref{eq:def of Ay_alpha}. By utilizing \eqref{eq:uniform convergence of S_inv}, \eqref{eq:uniform convergence of derivative of S}, \eqref{eq:limit of PPN} and the boundedness of $\sup_{\bm\eta\in\mathcal{D}_{\bm\eta}}\|{\partial^2 \mathbf{P}(\bm\eta)}/{\partial\eta_{k}\partial\eta_{l}}\|_{F}$, we obtain the uniform convergence of ${\partial^2 N\mathscr{F}_{\ty,\alpha}(\bm\theta_{0},\bm\eta)}/{\partial\eta_{k}\partial\eta_{l}}$ to $[\mathbf{A}_{\ty,\alpha}]_{k,l}$. The derivation of $\mathbf{B}_{\ty}$ is similar and thus omitted here.


Since $\mathscr{F}_{\Sy,\alpha}(\bm\eta)$ and $\mathscr{F}_{\GCV,\alpha}(\bm\eta)$ have the same dominant term as in Remark \ref{rmk:functions of ML}, i.e., $\mathscr{F}_{\ty,\alpha}(\hat{\bm\theta}^{\ML},\bm\eta)$ defined in Remark \ref{rmk:functions of ML}, the derivations for $\hat{\bm\eta}_{\GCV,\alpha}$ are similar to those of $\hat{\bm\eta}_{\Sy,\alpha}$ and omitted.

\subsection{Proof of Corollary \ref{corollary:sign of VarHPE of Sy GCV}}

For $\mathbf{P}(\eta)=\eta\mathbf{K}$ with fixed $\mathbf{K}\succ 0$, we have
\begin{align*}
\eta_{\ty\star,\alpha}=\frac{\bm\theta_{0}^{\top}\mathbf{K}^{-1}\mathbf{\Sigma}^{-1}\mathbf{K}^{-1}\bm\theta_{0}}{\alpha\Tr(\mathbf{\Sigma}^{-1}\mathbf{K}^{-1})}.
\end{align*}
It follows that for $\hat{\bm\theta}^{\TR}(\hat{\eta}_{\Sy,\alpha})$ or $\hat{\bm\theta}^{\TR}(\hat{\eta}_{\GCV,\alpha})$,
\begin{align}\label{eq:XVarHPE of SY and GCV}
&\Tr[\XVarHPE(\cdot)]=\nonumber\\
&\frac{4\alpha(\sigma^2)^2\Tr(\mathbf{\Sigma}^{-1}\mathbf{K}^{-1})}{(\bm\theta_{0}^{\top}\mathbf{K}^{-1} \mathbf{\Sigma}^{-1}\mathbf{K}^{-1}\bm\theta_{0})^2}\bm\theta_{0}^{\top}\mathbf{K}^{-1}\mathbf{\Sigma}^{-1}\mathbf{K}^{-1}\mathbf{\Sigma}^{-2}\mathbf{K}^{-1}\bm\theta_{0}.\qquad 
\end{align}
Even if $\mathbf{\Sigma}\succ 0$ and $\mathbf{K}\succ 0$, since $\bm\theta_{0}^{\top}\mathbf{K}^{-1}\mathbf{\Sigma}^{-1}\mathbf{K}^{-1}\mathbf{\Sigma}^{-2}\mathbf{K}^{-1}\bm\theta_{0}$ is not necessarily positive, it is hard to generally conclude if $\Tr[\XVarHPE(\cdot)]$ is positive or negative.
However, as long as $\mathbf{\Sigma}=\mathbf{I}_{n}$ or $\mathbf{K}=\mathbf{I}_{n}$, it is clear $\bm\theta_{0}^{\top}\mathbf{K}^{-1}\mathbf{\Sigma}^{-1}\mathbf{K}^{-1}\mathbf{\Sigma}^{-2}\mathbf{K}^{-1}\bm\theta_{0}$ is positive and thus $\Tr[\XVarHPE(\cdot)]>0$. 

For $\mathbf{\Sigma}=\diag\{\bm{s}\}$ and $\mathbf{P}(\bm\eta)=\diag\{\bm\eta\}$, we can prove that $[\bm\eta_{\tb\star,\alpha}]_{k}=[\bm\eta_{\ty\star,\alpha}]_{k}=g_{0,k}^2/\alpha$, and $\mathbf{b}_{\bm\eta,\star}^{\TR'}(\bm\eta)$, $\mathbf{A}_{\tb}^{-1}$, $\mathbf{B}_{\tb}$, $\mathbf{A}_{\ty}^{-1}$ and $\mathbf{B}_{\ty}$ are all diagonal with $[\mathbf{b}_{\bm\eta,\star}^{'}(\bm\eta)]_{k,k}={\sigma^2g_{0,k}}/{\eta_{k}^2s_{k}}$,
\begin{align*}
[\mathbf{A}_{\tb,\alpha}^{-1}]_{k,k}=&{\eta_{k}^3}/{(2g_{0,k}^2-\alpha\eta_{k})},\ [\mathbf{B}_{\tb}]_{k,k}=-{2g_{0,k}}/{\eta_{k}^2},\\
[\mathbf{A}_{\ty,\alpha}^{-1}]_{k,k}=&\frac{\eta_{k}^4s_{k}}{2(\sigma^2)^2(3g_{0,k}^2-2\alpha\eta_{k})},\
[\mathbf{B}_{\ty}]_{k,k}=-\frac{4(\sigma^2)^2g_{0,k}}{\eta_{k}^3s_{k}}.
\end{align*}
It follows that $\Tr[\XVarHPE(\cdot)]$ of $\hat{\bm\theta}^{\TR}(\hat{\bm\eta}_{\Sy,\alpha})$ or $\hat{\bm\theta}^{\TR}(\hat{\bm\eta}_{\GCV,\alpha})$ is $\sum_{k=1}^{p}[{4\alpha(\sigma^2)^2}/{g_{0,k}^2s_{k}^2}]>0$.

\section*{Acknowledgment} 

The authors thank Prof. Tianshi Chen for the stimulating discussion on the misalignment between the model parameters and the selected kernel. The authors also thank Dr. K\'evin Colin for his valuable feedback on this research.


\bibliographystyle{abbrv}
\bibliography{database} 

\vskip 0pt plus -1fil

\begin{IEEEbiography}[{\includegraphics[width=1in,height=1.25in,clip,keepaspectratio]{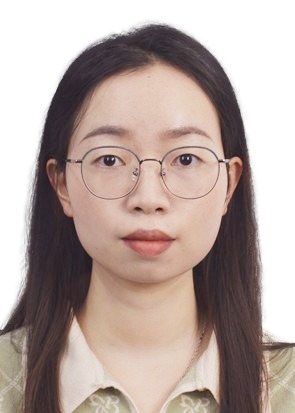}}] {Yue Ju} received the B.E. degree in automation from the Nanjing University of Science and Technology, Nanjing, China, in 2017 and the Ph.D. degree in computer and information en- gineering from the Chinese University of Hong Kong, Shenzhen, China, in 2022.
	
She is currently a postdoc at the KTH Royal Institute of Technology. She has been mainly working in the area of system identification.
\end{IEEEbiography}

\begin{IEEEbiography}[{\includegraphics[width=1in,height=1.25in,clip,keepaspectratio]{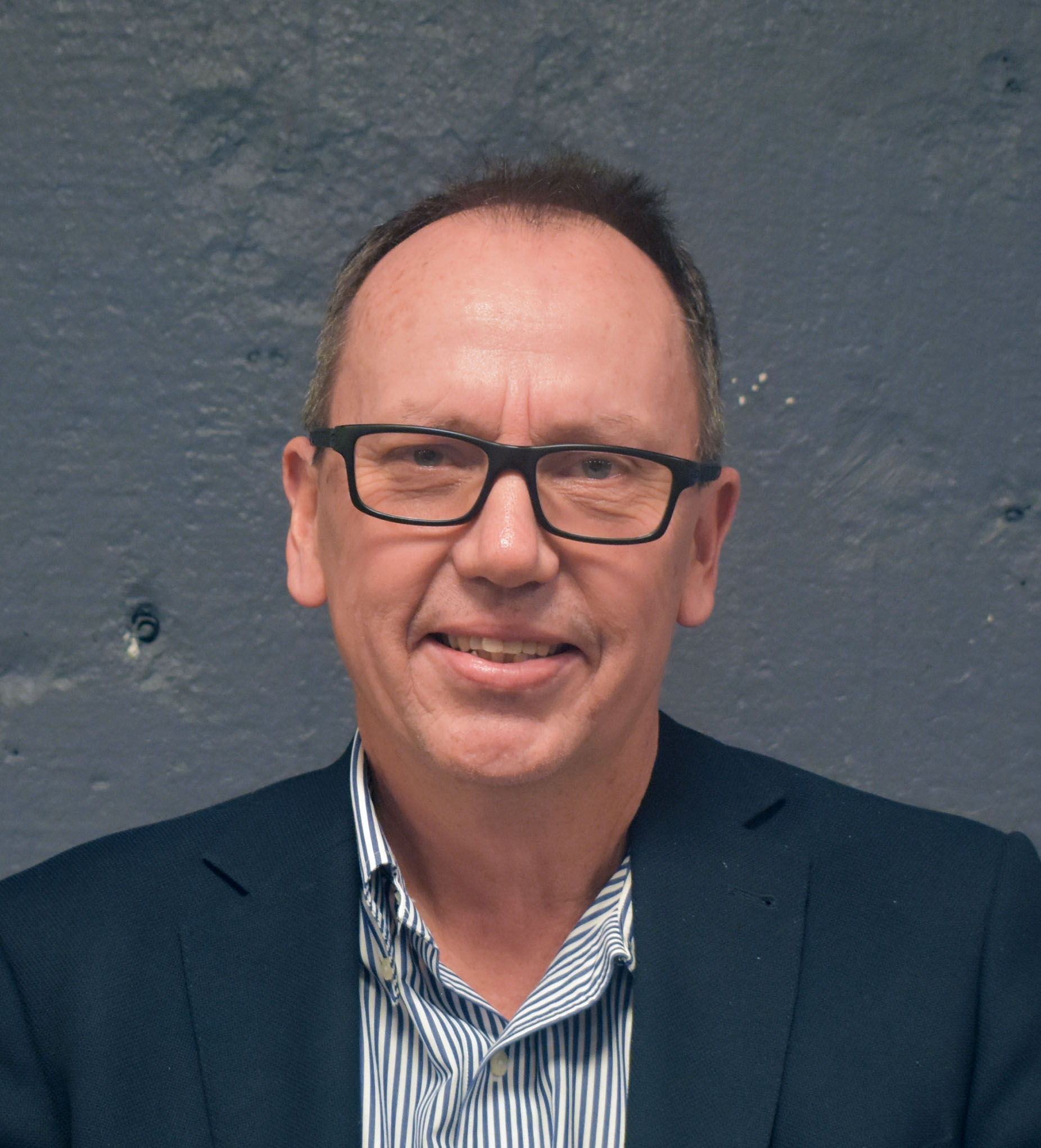}}]{Bo Wahlberg} has been the Professor of Chair in Automatic Control at KTH Royal Institute of Technology since 1991. He was elected an IEEE Fellow in 2007 for his contributions to system identification using orthonormal basis functions and a Fellow of IFAC in 2019 for his contributions to system identification and the development of orthonormal basis function models. Bo Wahlberg is since 2023 a Fellow of the Royal Swedish Academy of Engineering Sciences (IVA), Electrical Engineering Division.
 
Bo Wahlberg has received several awards, including the IEEE Transactions on Automation Science and Engineering Best New Application Paper Award in 2016. He is the author of over 250 scientific publications and has been the supervisor of more than 130 master’s students and 25 PhD students. His main research interest is in estimation and optimization in system identification, decision and control systems and signal processing with applications in process industry and transportation.
\end{IEEEbiography}

\vskip 0pt plus -1fil

\begin{IEEEbiography}[{\includegraphics[width=1in,height=1.25in,clip,keepaspectratio]{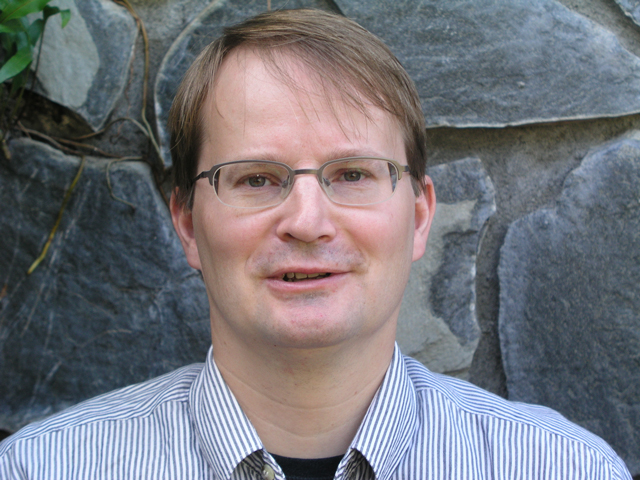}}] {H\r{a}kan Hjalmarsson} was born in 1962. He received the M.S. degree in Electrical Engineering in 1988, and the Licentiate degree and the Ph.D. degree in Automatic Control in 1990 and 1993, respectively, all from Linkˆping University,  Sweden.  He has held visiting research positions at California Institute of Technology, Louvain University and at the University of Newcastle, Australia.  He has served as an Associate Editor for Automatica (1996-2001), and IEEE Transactions on Automatic Control (2005-2007) and been Guest Editor for European Journal of Control and Control Engineering Practice. He is Professor at the Division of Decision and Control Systems, School of Electrical Engineering and Computer Science, KTH, Stockholm, Sweden and also affiliated with the Competence Centre for Advanced BioProduction by Continuous Processing, AdBIOPRO.  He is an IEEE Fellow and past Chair of the IFAC Coordinating Committee CC1 Systems and Signals. In 2001, he received the KTH award for outstanding contribution to undergraduate education. He was General Chair for the IFAC Symposium on System Identification held in 2018. His research interests include system identification, learning of dynamical systems for control, process modeling control and also estimation in communication networks. 
\end{IEEEbiography}

\end{document}